\documentclass[a4paper]{amsart}

\usepackage{enumerate,amssymb,amsfonts,latexsym,psfrag,ifthen,amsthm}
\usepackage{amscd,amsmath}
\usepackage{graphicx}

\newcommand{\Int}{\operatorname{Int}}
\newcommand{\id}{\operatorname{id}}
\newcommand{\bt}{\operatorname{bt}}
\newcommand{\snc}{\operatorname{snc}}
\newcommand{\MCG}{\operatorname{MCG}}
\newcommand{\Index}{\operatorname{index}}
\newcommand{\Fix}{\operatorname{Fix}}
\newcommand{\BT}{\text{BT}}

\newcommand{\Z}{\mathbb{Z}}
\newcommand{\R}{\mathbb{R}}
\newcommand{\Q}{\mathbb{Q}}
\newcommand{\N}{\mathbb{N}}

\newcommand{\cK}{\mathcal K}
\newcommand{\cD}{\mathcal D}

\newcommand{\wh}{\widehat}
\newcommand{\hw}{\hat{w}}
\newcommand{\ka}{\kappa}
\newcommand{\bh}{\bar{h}}
\newcommand{\I}{^{-1}}
\newcommand{\ol}[1]{\overline{#1}}

\newcommand{\ri}[1]{{#1}^\infty}
\newcommand{\li}[1]{{\vphantom{#1}}^\infty\! #1}
\newcommand{\rib}[1]{\ri{\left(#1\right)}}
\newcommand{\lib}[1]{\li{\left(#1\right)}}

\newcommand{\lp}[1]{\vphantom{#1}^+\! #1}

\newcommand{\rotint}{\rho i}

\newcommand{\opti}
{\raisebox{-.04cm}{$\stackrel{{\scriptscriptstyle 0}}{{\scriptscriptstyle 
1}}$}}

 \newcommand{\pichere}[2]
 {\begin{center}\includegraphics[width=#1\textwidth]{#2}\end{center}}

 \newcommand{\lab}[3]{\psfrag{#1}[#3]{$\scriptstyle{#2}$}
}

\newtheorem{thm}{Theorem}
\newtheorem{lem}[thm]{Lemma}
\newtheorem{cor}[thm]{Corollary}
\newtheorem{conj}[thm]{Conjecture}                            

\theoremstyle{definition}

\newtheorem{defn}[thm]{Definition}
\newtheorem{defns}[thm]{Definitions}
\newtheorem{example}[thm]{Example}

\theoremstyle{remark}

\newtheorem{remark}[thm]{Remark}
\newtheorem{remarks}[thm]{Remarks}

\begin{document}

\author{Andr\'e de Carvalho}
\address{Departamento de Matem\'atica Aplicada\\IME-USP\\Rua Do
  Mat\~ao 1010\\Cidade Universit\'aria\\05508-090 S\~ao Paulo
  SP\\Brazil}
\email{andre@ime.usp.br}

\author{Toby Hall}
\address{Department of Mathematical Sciences\\University of
  Liverpool\\Liverpool L69 7ZL\\UK}
\email{t.hall@liv.ac.uk}

\title[Decoration invariants]{Decoration invariants for horseshoe
  braids}

\begin{abstract}
The Decoration Conjecture describes the structure of the set of braid
types of Smale's horseshoe map ordered by forcing, providing
information about the order in which periodic orbits can appear when a
horseshoe is created. A proof of this conjecture is given for the
class of so-called {\em lone} decorations, and it is explained how to
calculate associated braid conjugacy invariants which provide
additional information about forcing for horseshoe braids.
\end{abstract}


\maketitle

\section{Introduction}

{\em Forcing relations} are a valuable tool in the dynamical study
of parameterized families of transformations: they provide information
about when the presence of certain dynamical features, such as the
existence of periodic orbits of a particular type, imply the presence
of other dynamical features.

For surface homeomorphisms, Smale's horseshoe is the paradigmatic map
with complicated dynamical behaviour, and understanding how it is
created in parameterized families of homeomorphisms is an important
problem. In this context, it is fruitful to study the forcing order on
the set of braid types of horseshoe periodic orbits~\cite{Boy}: this
partial order describes constraints on the order in which periodic
orbits can appear during the creation of a horseshoe.

Algorithmic implementations~\cite{BH,FM,Los} of Thurston's
classification of surface homeomorphisms provide a means of deciding
whether or not one given braid type forces another, but such an
approach doesn't provide any information about the global structure of
the forcing order on the set of all horseshoe braid types.

The {\em decoration conjecture}~\cite{frhbt} claims that the set of
horseshoe braid types is partitioned into families
$\cD^w=\{\beta^w_q\}$, each parameterized by a rational number~$q$,
which are totally ordered by the forcing order, in such a way that
$\beta^w_{q}$ forces $\beta^w_{q'}$ if and only if $q\le q'$. The
families are labelled by {\em decorations} $w$, which are finite words
in the symbols~$0$ and~$1$. Within families this trivializes the
computation of forcing: simply compare the rational parameters in the
usual order (or, equivalently, compare the symbolic representations of
the braids using the unimodal order).

Some special cases of this conjecture have been proved, and it is
supported by strong intuitive evidence from pruning theory: however, a
general proof has so far been elusive.  In this paper, the conjecture
is proved for a class of decorations called {\em lone}. There are many
lone decorations: two infinite families of them are described in
Sections~\ref{sec:star} and~\ref{sec:111} below. Collins~\cite{Piet}
states that 21 of the 63 decorations of length~5 or less are
lone\footnote{These are:
the empty decoration~$\cdot$, $0$, $1$, $00$, $11$, $000$, $111$, $101$,
$0000$, $0110$, $1111$, $1001$, $00000$, $01001$, $11001$, $10010$,
$10011$, $11011$, $11111$, $10101$, and $10001$.}.

The proofs of the main results presented here combine the pruning
techniques introduced in~\cite{Pruning} with unremovability arguments
of the type developed in~\cite{Ha}. Although pruning has provided an
inspiration for several previous results about forcing, this seems to
be the first time that it has been used effectively in proving such
results. 

\bigskip

The second goal of the paper is to present a new set of braid type
invariants.  Each totally ordered family $\cD^w = \{\beta^w_q\}$ gives
rise to an invariant (or, looking at it from the point of view of the
horseshoe braids themselves, a braid conjugacy invariant) $r^w$,
defined on the set of all horseshoe braid types~$\beta$ by
\[r^w(\beta)=\inf\{q\,:\,\beta^w_q\le\beta\}.\]
These invariants provide much additional information about the forcing
order: a periodic orbit of braid type~$\beta$ can only be created once
periodic orbits of braid types $\beta^w_q$ have been created for all
$q>r^w(\beta)$. Notice in this statement that, while $\beta^w_q$ is
restricted to be a horseshoe braid of lone decoration, $\beta$ can be
{\em any} horseshoe braid.

The techniques presented here make it possible, for lone
decorations, both to prove that the
family $\cD^w$ is totally ordered by forcing, and to calculate
the associated {\em decoration invariant} $r^w$.

\bigskip

The necessary background material on horseshoe braids and the
decoration conjecture is given in Section~\ref{sec:hsintro}, before
the main theorem and the algorithm for computing decoration invariants
are presented in Section~\ref{sec:results}. Some further background
material required for the proof, mostly concerning pruning and the
Asimov-Franks theorem, is given in Section~\ref{sec:af}.  The proof of
the main theorem is given in Section~\ref{sec:proof}. This is followed
by some examples and applications in Section~\ref{sec:examples}. The
applications include:
\begin{itemize}
\item a treatment of the so-called ``Star'' decorations
  (Section~\ref{sec:star}), which were introduced in~\cite{Stars},
  including the completion of the main theorem of that paper
  (Theorem~\ref{thm:starforce} here), thus providing a description of
  forcing {\em between} the families corresponding to different star
  decorations;
\item an example of how decoration invariants can be used to prove
  that certain other decorations are lone, and hence provide their own
  invariants (Theorem~\ref{thm:lone}); 
\item an example of how decoration invariants can be used to prove
  that certain horseshoe braids are of pseudo-Anosov type
  (Theorem~\ref{thm:pa}); and
\item a discussion of topological entropy bounds arising from
  decoration invariants (Section~\ref{sec:entropy}).
\end{itemize}

\section{Horseshoe braids: height and decoration}
\label{sec:hsintro}

This section contains the background material necessary to understand
the statement of the main theorem (Theorem~\ref{thm:main}). Although
complete, the treatment is rather terse: the
papers~\cite{frhbt,survey} are recommended for readers seeking a more
detailed account.

\subsection{Smale's horseshoe and the unimodal order}
\label{sec:smal-hors-unim}

\begin{figure}[htbp]
\lab{a}{H_0}{}
\lab{b}{H_1}{}
\lab{c}{S}{b}
\begin{center}
\pichere{0.4}{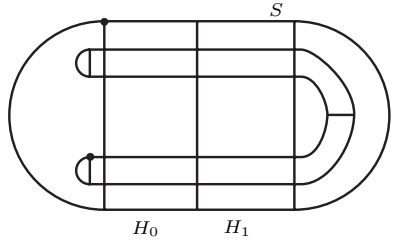}
\end{center}
\caption{Smale's horseshoe map}
\label{fig:hs}
\end{figure}

 In this paper, the standard model of Smale's horseshoe
 map~\cite{Smale} $F:D^2\to D^2$ depicted in Figure~\ref{fig:hs} is
 used. The set
\[\Lambda = \{x\in D^2\,:\,F^n(x)\in S\text{ for all }n\in\Z\}\]
(where~$S$ is the square depicted in Figure~\ref{fig:hs}) is a Cantor
set, and the {\em itinerary map} $k:\Lambda\to\{0,1\}^\Z$ defined by
\[k(x)_i=\left\{
\begin{array}{ll}
0 & \quad\text{ if }F^i(x)\in H_0\\
1 & \quad\text{ if }F^i(x)\in H_1
\end{array}
\right.\]
is a homeomorphism, conjugating $F|_{\Lambda}:\Lambda\to\Lambda$ to the {\em
  shift map} \mbox{$\sigma:\{0,1\}^\Z\to\{0,1\}^\Z$}.

The {\em unimodal order} $\preceq$ is a total order defined on
$\{0,1\}^\N$ as follows: if $s,t\in\{0,1\}^\N$, then $s\preceq t$ if
and only if either $s=t$, or the word $s_0s_1\ldots s_i$ contains an
even number of $1$s, where $i$ is least such that $s_i\not=t_i$. 

This order reflects the horizontal and vertical ordering of
points $x,y\in\Lambda$. Define the {\em horizontal and vertical
  coordinate functions} $h,v:\Lambda\to\{0,1\}^\N$ by
\mbox{$h(x)=k(x)_0k(x)_1k(x)_2\ldots$} and
$v(x)=k(x)_{-1}k(x)_{-2}k(x)_{-3}\ldots$. Then~$x$ lies to the left
of~$y$ if and only if $h(x)\prec h(y)$, and~$x$ lies below~$y$ if and
only if $v(x)\prec v(y)$. 

\subsection{Notation}

Much of the technical part of the paper is concerned with constructing
elements of $\{0,1\}^\N$ or $\{0,1\}^\Z$ from words $w\in\{0,1\}^j$
for various~$j$. The following notation will be used.

A {\em word} is an element of~$\bigcup_{j\ge 0}\{0,1\}^j$.

Let $u=u_0u_1\ldots u_{j-1}\in\{0,1\}^j$ be a word. Then $|u|=j$ denotes the
length of~$u$. $u$~is said to be {\em even (odd)} if it
contains an even (odd) number of $1$s. Let $\hat{u}=u_{j-1}\ldots
u_1u_0$ denote the reverse of~$u$, $\breve{u}=(1-u_0)u_1\ldots
u_{j-1}$ denote $u$ with the initial symbol changed, and
$\tilde{u}=u_0u_1\ldots (1-u_{j-1})$ denote $u$ with the final symbol
changed. When two of these accents are combined, they are applied
`bottom upwards': thus, for example,
$\widehat{\tilde{u}}=\breve{\hat{u}}$ is the word obtained by changing
the final symbol of~$u$ and then reversing the result. Denote by $u^+$
the word $u_0u_1\ldots u_{j-1}u_j\in\{0,1\}^{j+1}$, where $u_j$ is
chosen so that $u^+$ is even; and by~$\lp{u}$ the word
$u_{-1}u_0u_1\ldots u_{j-1}\in\{0,1\}^{j+1}$, where $u_{-1}$ is chosen
so that $\lp{u}$ is even.

$\overline{u}$ denotes the element $\ldots uuu\cdot uuu\ldots$ of
$\{0,1\}^\Z$ and $\ri{u}$ denotes the element $uuu\ldots$ of
$\{0,1\}^\N$. If $v,w,x$ are also words, then $\li{u}v\cdot w\ri{x}$
denotes the element $\ldots uuuv\cdot wxxx\ldots$ of $\{0,1\}^\Z$. If
$b=b_0b_1b_2\ldots$ and $f=f_0f_1f_2\ldots$ are elements of
$\{0,1\}^\N$, then $b\cdot f$ denotes the element $\ldots
b_2b_1b_0\cdot f_0f_1f_2\ldots$ of $\{0,1\}^\Z$. Similarly $bu\cdot
vf$ denotes the element \mbox{$\ldots b_2b_1b_0u\cdot vf_0f_1f_2\ldots$} of
$\{0,1\}^\Z$, $uf$ denotes the element $uf_0f_1f_2\ldots$ of
$\{0,1\}^\N$, and so on.

A {\em non-empty initial subword} of~$u\in\{0,1\}^j$ is a word
$u_0\ldots u_{i-1}\in\{0,1\}^i$ for some~$i$ with $0 < i\le
j$. Similarly, a {\em non-empty final subword} of $u$ is a word
$u_iu_{i+1}\ldots u_{j-1}$ for some~$i$ with $0\le i < j$.

\subsection{The horseshoe and its inverse}
\label{sec:hsinverse}
Recall (see~e.g.~\cite{conj}) that~$F$ is conjugate to its
inverse: $F^{-1}=\phi\circ F\circ \phi^{-1}$, where $\phi:D^2\to D^2$
is the (orientation-reversing) homeomorphism obtained by first
reflecting~$S$ in its horizonal centre line, and then rotating it
anticlockwise about its centre point through an
angle~$\pi/2$. 

The involution~$\phi$ restricts to an involution $\Lambda\to\Lambda$
which corresponds to reversing itineraries: if $k(x)=b\cdot f$ then
$k(\phi(x))=f\cdot b$.

\subsection{Height}

The {\em height function}~\cite{Ha} is a (not strictly) decreasing
function $\{0,1\}^\N\to[0,1/2]$ which is central to the results and
methods of this paper. In order to define it, is is necessary first to
introduce, for each rational $m/n\in(0,1/2]$, a word
  $c_{m/n}\in\{0,1\}^{n+1}$: these words will also play a central
  r\^ole throughout the paper.
\begin{defn}
\label{defn:cmn}
Let $m/n$ be a rational in $(0,1/2]$.  Let $L_{m/n}$ be the straight
  line in $\R^2$ from $(0,0)$ to $(n,m)$. For $0\le i\le n$, let
  $s_i=1$ if $L_{m/n}$ crosses some line $y=\text{integer}$ for
  $x\in(i-1,i+1)$, and $s_i=0$ otherwise. Then the word
  $c_{m/n}\in\{0,1\}^{n+1}$ is defined by $c_{m/n}=s_0s_1\ldots
  s_n$.  
\end{defn}

\begin{figure}[htbp]
\lab{0}{0}{}
\lab{1}{1}{}
\lab{2}{(0,0)}{r}
\lab{3}{(10,3)}{l}
\begin{center}
\pichere{0.5}{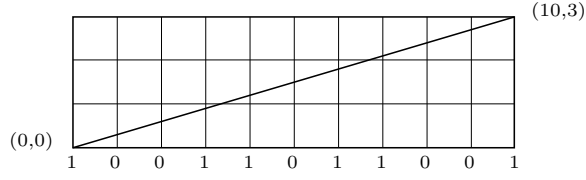}
\end{center}
\caption{$c_{3/10}=10011011001$}
\label{fig:cq}
\end{figure}

So, for example, $c_{3/10}=10011011001$ can be read off from
Figure~\ref{fig:cq}. The general form of these words can be seen in
Table~\ref{tab:prefix}, which shows~$c_{m/n}$ for all $m/n\in(0,1/2)$
with $1\le m\le 4$ and $3\le n\le 11$. Note that the words~$c_q$ are
clearly palindromic, and $c_{m/n}$ is of the form
$10^{\kappa_1}1^20^{\kappa_2}1^2\ldots 1^20^{\kappa_m}1$ for some
integers $\kappa_i\ge 0$ (an explicit formula for $\kappa_i$ can be
found in~\cite{Ha}, but is not needed here: note, however, that each
$\kappa_i$ is equal either to~$\kappa_1$ or to $\kappa_1-1$).  It can
also be seen easily from this description (see Lemma~2.7 of~\cite{Ha})
that if $m/n<m'/n'$, then $(c_{m'/n'}0)^\infty \prec
(c_{m/n}0)^\infty$.

\begin{table}[htbp]
\begin{center}
\begin{tabular}{|c||c|c|c|c|}
\hline
 & 1 & 2 & 3 & 4 \\
\hline \hline
3 & 1001 & & & \\
\hline
4 & 10001 & & & \\
\hline
5 & 100001 & 101101 & & \\
\hline
6 & 1000001 & & & \\
\hline
7 & 10000001 & 10011001 & 10111101 & \\
\hline
8 & 100000001 & & 101101101 & \\
\hline
9 & 1000000001 & 1000110001 & & 1011111101 \\
\hline
10 & 10000000001 & & 10011011001 & \\
\hline
11 & 100000000001 & 100001100001 & 100110011001 & 101101101101 \\
\hline
\end{tabular}
\end{center}
\caption{Examples of the words $c_{m/n}$ ($1\le m\le 4$, $3\le n\le
  11$)}
\label{tab:prefix}
\end{table}

\begin{defn}
\label{defn:height}
Let $c\in\{0,1\}^\N$. The {\em height} $q(c)$ of $c$ is the unique
element of~$[0,1/2]$ with the property that $c \prec (c_q0)^\infty$
for all rationals $q\in(0,q(c))$, and $(c_q0)^\infty \prec c$ for all
rationals $q\in(q(c),1/2]$.
\end{defn}

It is clear from the definition that $q:\{0,1\}^\N\to[0,1/2]$ is
decreasing. The next lemma (Theorem~3.2 of~\cite{Ha}) provides a
practical means of calculating $q(c)$ for all $c\in\{0,1\}^\N$ which
contain the subword $010$ (this is all that will be required in this
paper). 

\begin{lem}
\label{lem:heightcalc}
Let $c\in\{0,1\}^\N$, and suppose that $c$ contains the subword
$010$. Then $q(c)$ is rational, and can be calculated as
follows. First, if $10$ is not an initial subword of~$c$,
then~$q(c)=1/2$. If~$10$ is an initial subword, then write
\[c=10^{\kappa_1}1^{\mu_1}0^{\kappa_2}1^{\mu_2}\ldots,\]
where each $\kappa_i\ge 0$, each $\mu_i$ is either $1$ or $2$, and
$\mu_i=1$ only if $\kappa_{i+1}>0$ (thus $\kappa_i$ and $\mu_i$ are
uniquely determined by~$c$). For each $r\ge 1$, define
\[I_r(c)=\left(
\frac{r}{2r+\sum_{i=1}^r\kappa_i},
\frac{r}{2r-1+\sum_{i=1}^r\kappa_i}
\right],
\]
and let $s\ge 1$ be the least integer such that either $\mu_s=1$, or
$\bigcap_{r=1}^{s+1}I_r(c)=\emptyset$. Write $(x,y] =
  \bigcap_{r=1}^sI_r(c)\not=\emptyset$. Then
\[
q(c)=\left\{
\begin{array}{ll}
x & \quad\text{ if }\mu_s=2\text{ and } w\le x\text{ for all }w\in
I_{s+1}(c), \\ 
y & \quad \text{ if }\mu_s=1\text{, or }\mu_s=2\text{ and }w>y \text{
  for all }w\in I_{s+1}(c).
\end{array}
\right.
\]
\end{lem}

\begin{example}
Notice that the fact that $c$ contains the subword $010$ means that
$\mu_s=1$ for some~$s$, and hence the algorithm terminates. It may
also terminate before reaching the subword $010$. For example, let
$c=1011110011\ldots$. Then $\kappa_1=1$, $\mu_1=2$, $\kappa_2=0$,
$\mu_2=2$, $\kappa_3=2$, and $\mu_3=2$. This gives $I_1(c)=(1/3,
1/2]$, $I_2(c)=(2/5, 2/4]$, and $I_3(c)=(3/9, 3/8]$. Since $3/8<2/5$,
      the algorithm terminates and $q(c)=2/5$.
\end{example}

The following technical lemma (which follows immediately from Lemma~63
of~\cite{Stars}) will be needed:
\begin{lem}
\label{lem:starlem}
Let~$m/n$ be a rational in $(0,1/2]$, and let~$\kappa_i$ be integers
  such that 
\[c_{m/n}=10^{\kappa_1}1^20^{\kappa_2}1^2\ldots 1^20^{\kappa_m}1.\]
Let $1\le r\le m$, and let $f$ be any element of
$\{0,1\}^\N$. Then
\[q\left(10^{\kappa_r+1}1^20^{\kappa_{r+1}}1^2\ldots
1^20^{\kappa_m}1f\right) \le \frac{m}{n}.\]
\end{lem}

\subsection{Periodic orbits of the horseshoe}
Let~$P$ be a period~$n$ orbit of the horseshoe $F:D^2\to D^2$
(throughout this paper, ``period~$n$'' means {\em least}
period~$n$). Let $p$ be the rightmost point of~$P$: thus
$k(p)=\overline{c_P}$ for some word~$c_P$ of length~$n$, which is
called the {\em code} of the periodic orbit~$P$. Note that the choice
of~$p$ as the rightmost point of~$P$ means that
\[\sigma^i(\overline{c_P}) \prec \overline{c_P} \quad\text{ for } 1\le
i < n.\]

Recall~\cite{Boy} that the {\em braid type} $\bt(P\,;\,f)$ of a
period~$n$ orbit~$P$ of an orientation-preserving homeomorphism
$f:D^2\to D^2$ is a conjugacy class in the mapping class group
$\MCG(D_n)$ of the $n$-punctured disk~$D_n$, namely the conjugacy
class of the isotopy class of $h\I fh:D_n\to D_n$, where $h:D_n\to
D^2\setminus P$ is any orientation-preserving homeomorphism (if
$P\subseteq\partial D^2$, then first extend~$f$ over an exterior
collar). Braid types can thus be classified using the Thurston
classification~\cite{Thurston} as {\em finite order}, {\em reducible},
or {\em pseudo-Anosov}. The {\em forcing order} $\le$ on the set $\BT$
of all braid types is a partial order defined as follows: if
$\beta,\gamma\in\BT$, then $\beta\le\gamma$ if and only if every
orientation preserving homeomorphism $f:D^2\to D^2$ which has a
periodic orbit~$P$ with $\bt(P\,;\,f)=\gamma$ also has a periodic
orbit~$Q$ with $\bt(Q\,;\,f)=\beta$.

If~$P$ is a periodic orbit of the horseshoe, then the symbol~$P$
will often be used to denote the braid type $\bt(P\,;\,F)$ as well as
the periodic orbit itself. In particular, the notation $P\le Q$ is
used as a shorthand for $\bt(P\,;\,F)\le\bt(Q\,;\,F)$.

Most periodic orbits~$P$ of the horseshoe are {\em paired}, in the
sense that $\widetilde c_P$ is also the code of a horseshoe periodic
orbit~$\widetilde P$: in this case $\bt(P\,;\,F)=\bt(\widetilde
P\,;\,F)$. Thus, for example, the two periodic orbits~$P$
and~$\widetilde P$ with codes $c_P=10010$ and $c_{\tilde P}=10011$
have the same braid type, and it is common to write $c_P=1001\opti$,
reflecting the fact that the object of interest is the braid type
rather than the periodic orbit itself. The only orbits which are not
paired are those of even period $2k$, whose codes are of the form
$c_P=w\tilde w$ for some word~$w$ of length~$k$.

The {\em height} $q(P)\in(0,1/2]$ of a horseshoe periodic orbit~$P$ is
  defined to be $q(c_P^\infty)$. It is a braid type
  invariant~\cite{Ha}, so in particular $q(P)=q(\widetilde P)$ if~$P$
  is paired, and (taking the code which ends with~0), $q(P)$ is
  rational and can be calculated using the algorithm of
  Lemma~\ref{lem:heightcalc}. If~$P$ is not paired, then again
  $c_P^\infty$ contains the word $010$, so the algorithm terminates
  and $q(P)$ is rational.

Periodic orbits of the horseshoe can be classified as follows:
\begin{description}
\item[Orbits of finite order braid type] There are two fixed points,
  with codes $0$ and $1$ (and a fixed point outside of~$S$). There is
  one period two orbit, with code~$10$.

  For each rational $m/n\in (0,1/2)$, there is exactly one pair of
  period~$n$ orbits of finite order braid type whose rotation number
  about the fixed point of code~$1$ is $m/n$. The codes of these
  orbits are $d_{m/n}\opti$, where $d_{m/n}$ is the word consisting of
  the first $n-1$ symbols of $c_{m/n}$. These periodic orbits have
  height~$m/n$. There are no other periodic orbits of finite order
  braid type. (These results are due to Holmes and
  Williams~\cite{HoWi}.)
\item[NBT orbits] For each rational $q=m/n\in(0,1/2)$, the words
  $c_{m/n}\opti$ are the codes of a pair of period~$n+2$ orbits of
  pseudo-Anosov braid type, called {\em NBT orbits} and
  denoted~$P^\ast_q$~\cite{Ha}. These periodic orbits have
  height~$q$. There are no other horseshoe periodic orbits of these
  braid types. For $q=1/2$, the word $c_{1/2}1=1011$ is the code of a
  periodic orbit of reducible braid type, denoted~$P^\ast_{1/2}$.
\item[Orbits described by height and decoration] All other horseshoe
  periodic orbits~$P$ can be described by their height
  $q(P)\in(0,1/2]\cap\Q$ and their {\em decoration}, which is a word
    $w\in\{0,1\}^k$ for some $k\ge 0$: a periodic orbit of height~$q$
    has decoration~$w$ if and only if it has one of the four codes
    $c_q\opti w\opti$. Two periodic orbits of the same height and
    decoration have the same braid type~\cite{conj}: the
    notation~$P^w_q$ can therefore be used for any periodic orbit of
    height~$q$ and decoration~$w$.

    Only certain heights are compatible with a given
    decoration~$w$. Define the {\em scope}~$q_w$ of~$w$ to be the
    height of the horseshoe periodic orbit containing the point with
    itinerary $\overline{10w0}$, that is,
    \[q_w=\min_{0\le i\le
      k+2}q\left(\sigma^i(\rib{10w0})\right).\] The following result
    can be found in~\cite{frhbt}.
    \begin{lem}
      \label{lem:qw}
      For $q<q_w$ there are four periodic orbits of height~$q$ and
      decoration~$w$ (i.e. those with codes $c_q\opti w\opti$), while
      for $q>q_w$ there are no periodic orbits of height~$q$ and
      decoration~$w$. (For $q=q_w$, the four words $c_q\opti w\opti$
      may or may not be codes of periodic orbits of height~$q$.)
    \end{lem}
\end{description}

It is convenient for many purposes to consider~$\ast$ to be a
decoration as well (with scope $q_{\ast}=1/2$). Then every periodic
orbit~$P$ of the horseshoe is either of finite order braid type, or
can be described uniquely (up to braid type) by its
decoration~\mbox{$w\in\{\ast\}\cup\bigcup_{k\ge 0}\{0,1\}^k$} and its height
$q\le q_w$. (If~$q(P)=m/n$, then the period~$N$ of~$P$ is either $n$,
or greater than or equal to $n+2$. If $N=n$, then~$P$ is of finite
order braid type; if $N=n+2$, then~$P$ has decoration~$\ast$; and if
$N\ge n+3$, then~$P$ has decoration~$w$ of length $N-(n+3)$.)

The following lemma, which will be used several times, follows
immediately from the definition of~$q_w$ above, and the fact that
$q(\widehat{P})=q(P)$ for any horseshoe periodic orbit~$P$, where
$\widehat{P}$ is the periodic orbit containing the point of itinerary
$\overline{\widehat{c_P}}$ (Lemma~3.8 of~\cite{Ha}).
\begin{lem}
\label{lem:revq}
Let $w\in\bigcup_{k\ge 0}\{0,1\}^k$ be a decoration. Then $q_w=q_{\wh{w}}$.
\end{lem}

The results presented in this paper are related to the {\em Decoration
  Conjecture}~\cite{frhbt}. (Statement~c) below is not relevant in
this paper, and the definition of {\em topological train track type}
is therefore not given.)
\begin{conj}[Decoration Conjecture]
\label{conj:dec}
\mbox{}
\begin{enumerate}[a)]
\item There is an equivalence relation $\sim$ on the set of
  decorations, with the property that
  $\bt(P^w_q\,;\,F)=\bt(P^{w'}_{q'}\,;\,F)$ if and only if $q=q'$ and
  $w\sim w'$.
\item If $q\not=q_w$ then $P^w_q$ has pseudo-Anosov braid type.
\item For each decoration~$w$, the periodic orbits $P^w_q$ with
  $0<q<q_w$ all have the same topological train track type.
\item For each decoration~$w$, the set of braid types of periodic
  orbits of decoration~$w$ is totally ordered by forcing, with
  $P^w_q\le P^w_{q'}$ if and only if $q\ge q'$.
\item There is a partial order $\preceq$ on the set of equivalence
  classes of decorations, with the property that if $q>q'$ and
  $w\preceq w'$ then $P^w_q \le P^{w'}_{q'}$.
\end{enumerate}
\end{conj}

The main theorem of this paper, Theorem~\ref{thm:main} below,
concerns~d) and~e) of Conjecture~\ref{conj:dec}, for a certain class
of decorations:

\begin{defn}
A decoration $w\in\bigcup_{k\ge0}\{0,1\}^k$ is said to be {\em lone}
if for all \mbox{$q\in(0,q_w)\cap\Q$}, the four horseshoe periodic orbits of
height $q$ and decoration $w$ are the only horseshoe periodic orbits
of their braid type.
\end{defn}

Theorem~\ref{thm:main} states that Conjecture~\ref{conj:dec}~d)
holds for lone decorations~$w$, when restricted to those $P^w_q$ which
have pseudo-Anosov braid type (so if~b) is true, then so is~d) for
lone decorations). It also elucidates the partial order~$\preceq$
of~e) by providing a means of determining, for any lone decoration~$w$
and any horseshoe periodic orbit~$R$, which of the orbits $P^w_q$
are forced by~$R$. 

\begin{remark}
If Conjecture~\ref{conj:dec}~a) holds, then if there is a single
$q\in(0,q_w)\cap\Q$ for which the four orbits of height~$q$ and
decoration~$w$ are the only horseshoe periodic orbits of their braid
type, then the same is true for all~$q\in(0,q_w)\cap\Q$. In this case,
the lone decorations are precisely those which are alone in their
$\sim$-equivalence classes.
\end{remark}

The results stated in the next lemma can all be found in~\cite{Ha}:
\begin{lem}
\label{lem:heightfacts}
Let $P$ and $Q$ be periodic orbits of the horseshoe.
\begin{enumerate}[a)]
\item Height is a braid type invariant.
\item If $q<q(P)$ then $P^*_q\ge P$, while if $q>q(P)$ then
  $P^*_q\not\ge P$.
\item If $P\ge Q$ then $q(P)\le q(Q)$.
\end{enumerate}
\end{lem}

\section{Statement of results}
\label{sec:results}

Let~$w\in\bigcup_{k\ge 0}\{0,1\}^k$ be a lone decoration. Let
\[\cD^w=\{\bt(P^w_q\,;\,F)\,:\,q\in(0,q_w]\cap\Q\},\] the set of
  braid types of horseshoe periodic orbits of decoration~$w$. Let
  $\cK^w$ be the subset of $\cD^w$ consisting of pseudo-Anosov braid
  types: this contains $\bt(P^w_q\,;\,F)$ for a dense set of
  $q\in(0,q_w]\cap\Q$ by Lemma~\ref{lem:dense} below (and, if
    Conjecture~\ref{conj:dec}~b) holds, for all~$q$ except
    possibly~$q_w$). 

The main results of this paper are: 

\medskip\medskip

\begin{itemize}
\item that $\cK^w$ is totally ordered by
forcing, with $P^w_q\le P^w_{q'}$ if and only if $q\ge q'$; and 
\item that there exists a practical algorithm to determine, for any
  horseshoe periodic orbit~$R$, the number $r^w(R)\in(0,q_w]$ with the
    property that
\begin{eqnarray*}
R \ge P^w_q & \quad \text{ if } & q>r^w(R),\text{ and }\\
R \not \ge P^w_q  & \quad \text{ if } & q<r^w(R).
\end{eqnarray*}
In particular, for each lone decoration~$w$, $r^w$ is a braid type
invariant defined on the set of all horseshoe periodic orbits.
\end{itemize}

\medskip\medskip

The algorithm to compute $r^w(R)$ is complicated to state, although it
is easily implemented and is computationally light\footnote{The {\em
    horseshoe calculator} script at
  \texttt{http://www.maths.liv.ac.uk/cgi-bin/tobyhall/horseshoe}
  implements this algorithm.}. It will therefore be described
informally and illustrated by examples first; a formal description
which is more suitable for use in later proofs will then be given. The
purpose of the algorithm is to decide for which values of~$q$ points
of~$R$ are contained in certain disks bounded by segments of stable
and unstable manifold through points of the orbits~$P^w_q$: this
is what determines which of the~$P^w_q$ are forced by~$R$
(Theorem~\ref{thm:forcing}).

Let~$R$ be any horseshoe periodic orbit, with code~$c_R$: if~$R$ is
paired, then choose the code ending with~$0$ (so that the algorithm of
Lemma~\ref{lem:heightcalc} can be used to compute the various heights
below). Let~$w$ be any word. Then $r^w(R)$ is given by
\[r^w(R) = \min(\lambda^w(R), \max(\mu^w(R), \nu^w(R))),\]
where $\lambda^w(R)$, $\mu^w(R)$, and $\nu^w(R)$ are elements of
$(0,q_w]\cap\Q$ which will now be described.
\subsection*{$\mu^w(R)$}
Recall that $\lp{w}$ is the word obtained by prepending one symbol to
the front of~$w$, in such a way that $\lp{w}$ is even. Let~$v$ be a
non-empty even final subword of $\lp{w}$, and seek all occurrences of
the words $\breve{v}\opti 10$ in one period of $\overline{c_R}$ (that
is, all occurrences of either $\breve{v}010$ or of $\breve{v}110$:
recall that $\breve{v}$ is the word obtained from~$v$ by changing its
initial symbol). For each such occurrence, compute the height of the
forward sequence in $\overline{c_R}$ starting at the final symbols
$10$ in the occurrence. $\mu^w(R)$ is the minimum of such heights
taken over all such occurrences and all non-empty even final
subwords~$v$ of $\lp{w}$. If there are no such occurrences, or if the
minimum of the heights is greater than~$q_w$, then $\mu^w(R)=q_w$.

\begin{figure}[htbp]
\begin{center}
\pichere{0.5}{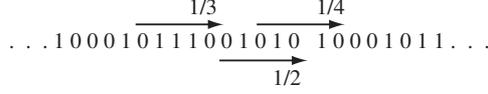}
\end{center}
\caption{Computing $\mu^w(R)$}
\label{fig:mu}
\end{figure}

\noindent\textbf{Example: } Let~$R$ have code~$c_R=100010111001010$,
and let~$w=1$. Then $\lp{w}=11$, which has only one non-empty even
final subword, namely $v=11$. Hence $\breve{v} = 01$, and occurrences
of $01010$ and $01110$ are sought. The three occurrences of such words
in one period of~$\overline{c_R}$ are shown in
Figure~\ref{fig:mu}. The corresponding forwards sequences have heights
$q(10010\ldots) = 1/3$, $q(1010\ldots) = 1/2$, and
$q(100010\ldots)=1/4$. Hence $\mu^w(R) = 1/4$.

\subsection*{$\nu^w(R)$}
The computation of $\nu^w(R)$ is similar: let~$v$ be a non-empty even
initial subword of $w^+$, and seek all occurrences of the words
$01\opti\tilde{v}$ in one period of $\overline{c_R}$. For each such
occurrence, compute the height of the backward sequence in
$\overline{c_R}$ starting at the initial symbols $01$ in the
occurrence. $\nu^w(R)$ is the minimum of such heights taken over
all such occurrences and all non-empty even initial subwords $v$ of
$w^+$. If there are no such occurrences, or if the minimum of the
heights is greater than~$q_w$, then $\nu^w(R)=q_w$.

\begin{figure}[htbp]
\begin{center}
\pichere{0.5}{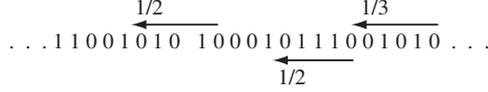}
\end{center}
\caption{Computing $\nu^w(R)$}
\label{fig:nu}
\end{figure}

\noindent\textbf{Example: } Continuing with the above example,
$w^+=11$ which has only one non-empty even initial subword, namely
$v=11$. Hence $\tilde{v} = 10$, and occurrences of $01010$ and $01110$
are sought. This gives $\nu^w(R)=1/3$ (see Figure~\ref{fig:nu}).

\subsection*{$\lambda^w(R)$}
Seek all occurrences of the four words $01\opti w\opti 10$ in one
period of $\overline{c_R}$. For each such occurrence, compute
$\max(q(b), q(f))$, where $b$ is the backward sequence starting at the
initial symbols $01$ in the occurrence, and $f$ is the forward sequence
starting at the final symbols $10$. Then $\lambda^w(R)$ is the minimum
value of this quantity taken over all such occurrences. If there are no
such occurrences, or if the minimum is greater than~$q_w$, then
$\lambda^w(R)=q_w$. 

\begin{figure}[htbp]
\begin{center}
\pichere{0.5}{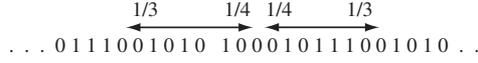}
\end{center}
\caption{Computing $\lambda^w(R)$}
\label{fig:lambda}
\end{figure}

\noindent\textbf{Example: } Continuing with the above example,
occurrences of $0101010$, $0111010$, $0101110$, and $0111110$ are
sought. There are two such occurrences, one with $q(b)=1/3$ and
$q(f)=1/4$, and the other with $q(b)=1/4$ and $q(f)=1/3$ (see
Figure~\ref{fig:lambda}). Hence $\max(q(b),q(f))=1/3$ for both
occurrences, and hence $\lambda^w(R)=1/3$.

Thus $r^w(R)=\min(\lambda^w(R), \max(\mu^w(R), \nu^w(R))) =
\min(1/3, \max(1/4, 1/3))= 1/3$.

\begin{example}
Table~\ref{tab:decinvs} lists all of the horseshoe periodic orbits of
period~8, together with the values of their decoration invariants for
the lone decorations of length~3 or less (and for the
decoration~$\ast$: see Remark~\ref{remark:NBTdecinv}). Period~8 orbits
with the same height and decoration (and hence the same braid type)
are grouped in the same row. Those orbits for which there is not a
choice of 4 different codes are either of finite order type
($1011011\opti$ and $1000000\opti$), NBT orbits ($1000001\opti$), or
orbits whose height is equal to the scope of their decoration
($10111010$, with height 1/2 and decoration $101$; $1011111\opti$ with
height~1/2 and decoration~$111$; $1001101\opti$, with height~1/3 and
decoration~01; and $10001\opti0(\opti)$, with height~1/4 and
decoration~0 --- here there are three orbits with this height and
decoration, namely those with codes $10001001$, $10001100$, and
$10001101$).

This table gives exact information about which of the (infinitely
many) orbits with these decorations are forced by which period~8
orbits.

Thus, for example, 7 of the orbit types have $r^{111}=1/2=q_{111}$, and
hence force no orbits of decoration $111$. The other 4 orbit types all
force some of the orbits with decoration $111$: for example, the orbits
of code~$10001\opti1\opti$ force $P^{111}_q$ for all $q>1/4$, and do not
force $P^{111}_q$ for all $q<1/4$.
\begin{table}[htbp]
\begin{center}
\begin{tabular}{|l|c|c|c|c|c|c|c|c|c|}
\hline
Decoration & $\ast$ & $\cdot$ & $0$ & $1$ & $00$ & $11$ &
$000$ & $101$ & $111$ \\
\hline
Scope & 1/2 & 1/3 & 1/4 & 1/2 & 1/5 & 2/5 & 1/6 & 1/2 & 1/2 \\
\hline
$10111010$ & 1/2 & 1/3 & 1/4 & 1/2 & 1/5 & 2/5 & 1/6 & 1/2 & 1/2 \\
$1011111\opti$ & 1/2 & 1/3 & 1/4 & 1/2 & 1/5 & 2/5 & 1/6 & 1/2 & 1/2 \\
$1011011\opti$& 1/2 & 1/3 & 1/4 & 1/2 & 1/5 & 2/5 & 1/6 & 1/2 & 1/2 \\
$1001\opti11\opti$ & 1/2 & 1/3 & 1/4 & 1/2 & 1/5 & 1/3 & 1/6 & 1/2 & 1/2 \\
$1001\opti10\opti$ & 1/3 & 1/3 & 1/4 & 1/3 & 1/5 & 1/3 & 1/6 & 1/3 & 1/3 \\
$1001101\opti$ & 1/3 & 1/3 & 1/4 & 1/3 & 1/5 & 1/3 & 1/6 & 1/3 & 1/3 \\
$10001\opti0(\opti)$ & 1/2 & 1/3 & 1/4 & 1/2 & 1/5 & 2/5 & 1/6 & 1/2 & 1/2 \\
$10001\opti1\opti$ & 1/2 & 1/4 & 1/4 & 1/4 & 1/5 & 1/4 & 1/6 & 1/2 & 1/4 \\
$100001\opti\opti$ & 1/2 & 1/5 & 1/5 & 1/2 & 1/5 & 2/5 & 1/6 & 1/2 & 1/2 \\
$1000001\opti$ & 1/6 & 1/6 & 1/6 & 1/6 & 1/6 & 1/6 & 1/6 & 1/6 & 1/6 \\
$1000000\opti$ & 1/2 & 1/3 & 1/4 & 1/2 & 1/5 & 2/5 & 1/6 & 1/2 & 1/2 \\
\hline
\end{tabular}
\end{center}
\caption{Examples of some decoration invariants for period 8 orbits}
\label{tab:decinvs}
\end{table}
 
On the whole, the orbits in this table force relatively few of the
orbits of the given decorations. This is because the period is
low. The following ``weak universality'' result says that, for a fixed
decoration~$w$, most horseshoe periodic orbits~$R$ have $r^w(R)<q^w$:
indeed, most have $r^w(R)$ very small, and hence force almost all
orbits of decoration~$w$.

\begin{thm}
\label{thm:universal}
Let~$w$ be a decoration and $q\in(0,q_w)$. Let $p_n$ denote the
proportion of period~$n$ horseshoe orbits~$R$ with $r^w(R)<q$. Then
$p_n\to 1$ as $n\to\infty$.
\end{thm}
\begin{proof}
Pick $k$ with $1/k<q$. Then any~$R$ whose code includes the word
$0^k10w010^k$ has
\[r^w(R) \le \lambda^w(R) < q(10^k\ldots) < 1/k < q.\]
\end{proof}

As a further illustration of the use of these invariants,
let~$R=P_{1/6}^{10}$ be the periodic orbit of
code~$c_R=10000011100$. Computing $r^w(R)$ for the same lone
decorations~$w$ as in Table~\ref{tab:decinvs} gives $r^{\ast}(R)=1/3$,
$r^{\cdot}(R)=1/3$, $r^{0}(R)=1/6$, $r^{1}(R)=1/3$, $r^{00}(R)=1/6$,
$r^{11}(R)=1/3$, $r^{000}(R)=1/6$, $r^{101}(R)=1/3$, and
$r^{111}(R)=1/3$. The periodic orbits of these decorations are
depicted in Figure~\ref{fig:forceshow}: for each decoration~$w$, there
are periodic orbits $P_q^w$ for each rational~$q$ in $(0,q_w)$
(represented by points on the vertical lines), but not for rationals
$q>q_w$. The decoration invariants show that all orbits on the thicker
parts of the lines are forced by~$R$, but that none of
the orbits on the thinner parts of the lines are forced by~$R$ (the
theorem doesn't state whether or not the orbits represented by points
at the transition from thin to thick are forced by~$R$). Notice
that since $q(R)=1/6$, $R$ cannot force any orbit $P_q^w$ with $q<1/6$
by Lemma~\ref{lem:heightfacts}~c). The decoration $10$ of $R$ itself
is {\em not} lone, and is not required to be by Theorem~\ref{thm:main}.

\begin{figure}[htbp]
\lab{q}{q}{b}
\lab{r}{R}{bl}
\begin{center}
\pichere{0.75}{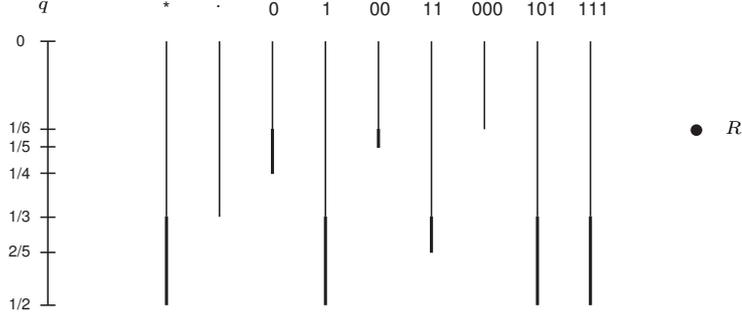}
\end{center}
\caption{Some orbits forced, and not forced, by the orbit of code
  10000011100}
\label{fig:forceshow}
\end{figure}
\end{example}

A formal statement of the algorithm for computing $r^w(R)$ will now be
given. Some preliminary definitions are useful.

\begin{defns}
\label{defn:forwardbackward}
Let $v$ and $c$ be words, with $|c|\ge 1$. For each $i$ with
$0\le i < |c|$, let
\[
\overrightarrow{r}(\overline{c},v,i)=
\left\{
\begin{array}{ll}
q(f) & \text{ if }\sigma^i(\overline{c})=bv\cdot f \text{
  for some }b,f\in\{0,1\}^\N\\
\frac{1}{2} & \text{ otherwise.}
\end{array}
\right.
\]
Similarly, define $\overleftarrow{r}(\overline{c},v,i)$ by replacing
$q(f)$ with $q(b)$ in the above, and
$\overleftrightarrow{r}(\overline{c},v,i)$ by replacing $q(f)$ with
$\max(q(b),q(f))$. 

Then, if~$V$ is a finite non-empty set of words, let
\[
\overrightarrow{r}(\overline{c}, V) = 
\min_{v\in V,\,\, 0\le i < |c|} \overrightarrow{r}(\overline{c},v,i).
\]
 $\overleftarrow{r}(\overline{c}, V)$ and
$\overleftrightarrow{r}(\overline{c}, V)$ are defined similarly using
$\overleftarrow{r}(\overline{c},v,i)$ and
$\overleftrightarrow{r}(\overline{c},v,i)$ in place of
$\overrightarrow{r}(\overline{c},v,i)$.

Finally, if~$w$ is a word and~$R$ is a horseshoe periodic orbit, define
\[\mu^w(R)=\min\left(
q_w, 
\overrightarrow{r}^w(\overline{c_R}, V)
\right),
\]
where
\[
V = \{
\breve v\opti\,:\,\text{$v$ is a non-empty even final subword of }\lp{w}
\}.
\] 

Similarly, define
\[\nu^w(R)=\min\left(
q_w, 
\overleftarrow{r}^w(\overline{c_R}, V)
\right),
\]
where
\[
V = \{
\opti\tilde v\,:\,\text{$v$ is a non-empty even initial subword of }w^+
\},
\]
and
\[\lambda^w(R)=\min\left(
q_w, 
\overleftrightarrow{r}^w(\overline{c_R}, V)
\right),
\]
where
\[
V = \{
0w0, 0w1, 1w0, 1w1
\}.
\]
\end{defns}

\begin{remark}
\label{remark:rreverse}
It is clear from the definitions that
$\lambda^{\hat{w}}(\widehat{R})=\lambda^w(R)$,
$\mu^{\hat{w}}(\widehat{R})=\nu(R)$, and
$\nu^{\hat{w}}(\widehat{R})=\mu(R)$
\end{remark}

 The pieces are now all in place to describe the invariant $r^w$:
 \begin{defn}
 \label{defn:rw}
 Let $w\in\bigcup_{k\ge 0}\{0,1\}^k$ be a decoration, and~$R$ be a
 horseshoe periodic orbit. Then the {\em $w$-depth of~$R$},
 $r^w(R)\in(0,q_w]\cap\Q$ is defined by
 \[r^w(R)=\min\left(
 \lambda^w(R), \max(\mu^w(R), \nu^w(R))
 \right).
 \]
\end{defn}

A number of the results presented below require that the braid type of
$P_q^w$ be pseudo-Anosov, hence the following definition:

\begin{defn}
Let $w$ be a decoration. Define 
\[\Q^w=\{q\in(0,q_w)\cap\Q\,:\,P_q^w\text{ has pseudo-Anosov braid
type}\}.\]
\end{defn}

Conjecture~\ref{conj:dec}~b) states that $\Q^w=(0,q_w)\cap\Q$ for all
decorations $w$, and this can be proved in a number of special
cases. However, in this paper all that is needed is the following
straightforward lemma: 

\begin{lem}
\label{lem:dense}
For any decoration $w$, $\Q^w$ is dense in $(0,q_w)$.
\end{lem}

\begin{proof}
 Let $w$ be a decoration of length $k\ge0$, and let
 $q=m/n\in(0,q_w)$. Since $P_q^w$ does not have finite order braid
 type, it must have pseudo-Anosov braid type if its period $n+k+3$
 is prime~\cite{Boy}. Hence
\[\Q^w\supseteq A^w=\{m/n\in(0,q_w)\,:\,(m,n)=1 \text{ and } 
n+k+3\text{ is prime}\}.\] However, $A^w$ is dense in
$(0,q_w)$. Indeed, for any $N\in\N$, the set \[\Q_N=\{m/n\,:\,(m,n)=1
\text{ and } n+N\text{ is prime}\}\] is dense in~$\Q$. For let
$r/s\in\Q$. For each $K\ge 1$, let $r_K=K(N-1)r+1$ and
\mbox{$s_K=K(N-1)s+1$}, so that $r_K/s_K\to r/s$ as $K\to\infty$, and it
suffices to show that each $r_K/s_K$ can be approximated arbitrarily
closely by elements of $\Q_N$. Since $r_Ks_K$ is coprime to $N-1$,
Dirichlet's theorem gives a strictly increasing sequence of integers $a_K^i$
such that $a_K^ir_Ks_K + (N-1)$ is prime for all~$i$. Let
$m_K^i=a_K^ir_K^2$ and $n_K^i = a_K^ir_Ks_K-1$. Then $m_K^i/n_K^i\to
r_K/s_K$ as $i\to\infty$, and $m_K^i/n_K^i\in\Q_N$ for all~$i$ as
required.
\end{proof}

The main theorem of this paper can now be stated:
\begin{thm}
\label{thm:main}
Let $w$ be a lone decoration. Then 
\begin{enumerate}[a)]
\item The braid types of the orbits $P_q^w$ with $q\in\Q^w$ are
  totally ordered by forcing, with \mbox{$P_q^w\ge P_{q'}^w$} if and
  only if $q\le q'$.
\item For any horseshoe periodic orbit $R$,
\[r^w(R)=\sup\{q\in\Q^w\,:\,R\not\ge P_q^w\},\]
and hence $r^w$ is a braid type invariant.
\end{enumerate}
\end{thm}

\begin{remark}
\label{remark:NBTdecinv}
There is a corresponding result for the family
$\{P^\ast_q\,:\,q\in(0,1/2]\cap\Q\}$ of NBT orbits~\cite{Ha}. This
  family is totally ordered by forcing (with $P_q^\ast \ge
  P_{q'}^\ast$ if and only if $q\le q'$), and the corresponding braid
  type invariant $r^\ast$ can be calculated by
\[r^\ast(R)=\overleftrightarrow{r}\left(\overline{c_R}, \{0,1\}\right).\]
\end{remark}

\begin{remark}
$r^w(R)$ gives information about which periodic orbits of
  decoration~$w$ are forced by~$R$. Invariants $q^w(R)$ which
  give information about which periodic orbits of decoration~$w$
  force~$R$ can be defined similarly:
\[q^w(R) = \sup\{q\in\Q^w\,:\, P^w_q \ge R\}\]
(or $q^w(R)=q_w$ if $P^w_q\not\ge R$ for all $q$). However, the
authors know of no means of computing these invariants, except for the
NBT decoration $w=\ast$, for which $q^\ast(R)=q(R)$ by
Lemma~\ref{lem:heightfacts}~b). 
\end{remark}

\section{Linking numbers, the Asimov-Franks theorem, and pruning}
\label{sec:af}

The main technique used in this paper to show that the braid type of
one periodic orbit~$P$ forces the braid type of a second periodic
orbit~$Q$ is to apply the Asimov-Franks theorem~\cite{AF,Unremovable}
to show that~$Q$ is {\em unremovable} in $D^2\setminus P$. In order to
show that the conditions of the Asimov-Franks theorem hold, linking
number considerations will be used. The necessary definitions and
results are presented in
Sections~\ref{sec:linking}\,--\,\ref{sec:method}.

On the other hand, the technique used to show that the braid type of
one periodic orbit does {\em not} force the braid type of a second is
to show that all orbits of the second braid type can be {\em pruned
  away} by an isotopy supported in the complement of the first
orbit. The pruning theorem used to do this~\cite{Pruning} is stated in
Section~\ref{sec:pruning}. 

\subsection{Linking numbers}
\label{sec:linking}
Let $f:D^2\to D^2$ be an orientation-preserving homeomorphism, and
$\{f_t\}:\id\simeq f$ be an isotopy (which will be referred to as a
{\em suspension} of~$f$). Recall that, given two distinct periodic
orbits~$P$ and~$Q$ of~$f$, the {\em linking number} $L(Q,P)$ of~$Q$
about~$P$ with respect to the suspension $\{f_t\}$ is an integer which
determines the homology class of the suspension
$\pi(\bigcup_{s\in[0,1]}(f_s(Q)\times\{s\}))$ of~$Q$ in the suspension
manifold
\[
\left[
\left(D^2\times[0,1]\right)\setminus\bigcup_{s\in[0,1]}(f_s(P)\times\{s\})
\right]
\biggl/(x,0)\sim(x,1).
\]

If the suspension $\{f_t\}$ is changed, then the linking number
$L(Q,P)$ changes by a multiple of the period~$n$ of~$Q$: thus the
linking number is well-defined modulo~$n$.

The following two straightforward lemmas can be found
in~\cite{thesis}, where they appear
as Lemma~1.23 and Corollary~1.25 respectively. The first can be used
to show that two periodic orbits have distinct linking numbers about a
third periodic orbit~$P$: intuitively, this means that they cannot
move together and annihilate under an isotopy relative to~$P$.

\begin{lem}
\label{lem:linkcalc}
Let $f\colon D^2\to D^2$ be an orientation preserving homeomorphism
which has periodic points~$q_0$ and~$q_1$ of least period~$n$ lying on
orbits $Q_0$ and $Q_1$, and let $\alpha\colon [0,1]\to D^2$ be a path
from~$q_0$ to~$q_1$. Let $\gamma$ be the closed curve
$\alpha\cdot(f^n\circ\alpha)\I$. Suppose that~$P$ is a periodic orbit
of~$f$ with no points lying on~$\alpha$, and that~$\gamma$ has winding
number $w_\gamma(P)$ about~$P$. Then
\[L(Q_0,P)=L(Q_1,P)+w_\gamma(P),\]
where the linking numbers are calculated with respect to any fixed
suspension of~$f$.
\end{lem}

The second lemma can be used to exclude the possibility of a given
periodic orbit collapsing onto an orbit of lower period under isotopy.

\begin{lem}
\label{lem:nocollapse}
Let $f_i:D^2\to D^2$ be a sequence of orientation-preserving
homeomorphisms, converging in the $C^0$ topology to a homeomorphism
$f:D^2\to D^2$. Let $p_i\to p$ and $q_i\to q$ be sequences in~$D^2$
with the property that each $p_i$ lies on a period~$m$ orbit~$P_i$
of~$f_i$ and each~$q_i$ lies on a period~$n$ orbit~$Q_i$
of~$f_i$. If~$q$ is a period~$n/l$ point of~$f$ which does not lie on
the $f$-orbit of~$p$, then $L(Q_i,P_i)$ is a multiple of~$l$ for
sufficiently large~$i$ (with respect to any suspension of~$f_i$).
\end{lem}

\subsection{The Asimov-Franks theorem}

The Asimov-Franks theorem~\cite{AF} gives conditions under which
periodic orbits of homeomorphisms $f:M\to M$ persist (or are {\em
  unremovable}) under arbitrary isotopy of~$f$. The version given here
is from~\cite{Unremovable}, and is restricted to the case of interest in
this paper (unremovability of single periodic orbits of
orientation-preserving homeomorphisms of $D^2$ under isotopy relative
to some other periodic orbit).

\begin{defns}
Let $f\colon D^2\to D^2$ and $g\colon D^2\to D^2$ be
orientation-preserving homeomorphisms having period $n$ points $p$ and
$q$ respectively; and let $R$ be a periodic orbit of~$f$ in
$\Int(D^2)$. Then $(p\,;\,f)$ and $(q\,;\,g)$ are {\em connected by
  isotopy rel. $R$} (denoted $(p\,;\,f)\sim(q\,;\,g)$) if
there exists an isotopy $\{f_t\}\colon f\simeq g$ relative to $R$ and
a path $\alpha$ in $D^2$ from $p$ to $q$, such that $\alpha(t)$ is a
period~$n$ point of $f_t$ for all $t$. 

$(p\,;\,f)$ is said to be {\em unremovable in $(D^2,R)$} if every
homeomorphism $g:D^2\to D^2$ which is isotopic to $f$ rel. $R$ has a
period~$n$ point~$q$ with $(p\,;\,f)\sim(q\,;\,g)$. 
\end{defns}

Notice that being connected by isotopy rel.~$R$ is clearly an
equivalence relation on the set of all pairs $(p\,;\,f)$,
where~$f:D^2\to D^2$ is an orientation-preserving homeomorphism
and~$p$ is a periodic point of~$f$.

The Asimov-Franks theorem gives three conditions which together ensure
the unremovability of~$(p\,;\,f)$. The first condition prevents the
orbit of~$p$ from collapsing onto an orbit of lower period under
isotopy.

\begin{defn}
$(p\,;\,f)$ is said to be {\em uncollapsible} if, given any sequence
  of homeomorphisms $g_j:D^2\to D^2$ which converge in the $C^0$
  topology to $g:D^2\to D^2$, and a sequence $q_j \to q$ in $D^2$ such
  that $q_j$ is a period~$n$ point of~$g_j$ with
  $(q_j\,;\,g_j)\sim(p\,;\,f)$ for all~$j$, then $q$ is a period~$n$
  point of~$g$.
\end{defn}

The second condition prevents the orbit of~$p$ from falling into the
periodic orbit~$R$. 

\begin{defn}
$(p\,;\,f)$ is said to be {\em separated from~$R$} if
  \mbox{$(p\,;\,f)\sim(q\,;\,g)$} implies $q\not\in R$.
\end{defn}

The final condition is that the total fixed point index of periodic
points which could interact with~$p$ under isotopy is non-zero. Recall
(see for example~\cite{Jiang}) that if $f:X\to X$ is a continuous
self-map of a compact manifold, then the fixed point index
$\Index(S,f)$ of a subset $S$ of $\Fix(f)$ can be defined,
generalising the familiar notion of the index of an isolated fixed
point, provided that~$S$ is compact and is open in~$\Fix(f)$.

\begin{defn}
Let~$p$ be a period~$n$ point of~$f$.  The {\em strong
  Nielsen class} $\snc(p\,;\,f)$ of $(p\,;\,f)$ is the set
of all period~$n$ points~$q$ of~$f$ with $(p\,;\,f)\sim(q\,;\,f)$. 
\end{defn}

If $(p\,;\,f)$ is uncollapsible and separated from~$R$, then
$\snc(p\,;\,f)$ is compact and open in the set of fixed points
of~$f^n$, so the following definitions can be made:
\begin{defns}
Let $(p\,;\,f)$ be uncollapsible and separated from~$R$.  The {\em
  index} of $(p\,;\,f)$ is defined by
$I(p\,;\,f)=\Index(\snc(p\,;\,f),f^n)$.  $(p\,;\,f)$ is said to be
{\em essential} if $I(p\,;\,f)\not=0$.
\end{defns}

\begin{thm}[Asimov-Franks]
\label{thm:AF}
If $(p\,;\,f)$ is uncollapsible, separated from~$R$, and essential,
then it is unremovable.
\end{thm}

The following result can be found
in~\cite{Unremovable}: it says that relevant topological information
is preserved under connection by isotopy.
\begin{lem}
\label{lem:icpreserve}
Let $(p\,;\,f)\sim(q\,;\,g)$, and let the orbits of~$p$ and~$q$ be
denoted~$P$ and~$Q$ respectively. Then
\begin{enumerate}[a)]
\item $\bt(P\,;\,f)=\bt(Q\,;\,g)$.
\item If $P\not=R$ and $Q\not=R$, then for any suspension of~$f$ there
  is a suspension of~$g$ such that $L(P,R)=L(Q,R)$ with respect to the
  given suspensions.
\end{enumerate}
\end{lem}

\begin{remark}
\label{rem:separated}
By Lemma~\ref{lem:icpreserve}~a), $(p\,;\,f)$ is
necessarily separated from~$R$ if
\mbox{$\bt(P\,;\,f)\not=\bt(R\,;\,f)$} (and in particular if~$P$
and~$R$ have different periods).
\end{remark}

The following trivial result will also be useful:
\begin{lem}
\label{lem:iciterate}
Let~$p$ and~$q$ be period~$n$ points of $f$ and $g$. Then
\mbox{$(p\,;\,f)\sim(q\,;\,g)$} if and only if
$(f(p)\,;\,f)\sim(g(q)\,;\,g)$.
\end{lem}

\subsection{A method for showing that two periodic points are
  connected by isotopy}
\label{sec:method}

\begin{defn}
\label{defn:define}
Let $f\colon D^2\to D^2$ be an orientation-preserving homeomorphism
with distinct period~$n$ points~$p$ and~$q$. Suppose that
$\alpha\colon[0,1]\to D^2$ is an arc from~$p$ to~$q$ with the property
that 
\begin{enumerate}[a)]
\item there exist $a,b\in[0,1]$ such that 
\begin{enumerate}[i)]
\item $\alpha([a,1])\cup
f^n(\alpha([b,1]))$ is a simple closed curve bounding a (closed) disk
$\Delta$, and 
\item $\alpha([0,a])=f^n(\alpha([0,b]))$, 
\end{enumerate}
\item $\alpha([0,1])$ is disjoint from $f^i(\alpha([0,1]))$ for $1\le
i<n$, and 
\item the orbits~$P$ and~$Q$ of~$p$ and~$q$ are disjoint from
$\Int\,\Delta$.
\end{enumerate}
Then $\alpha$ is said to {\em define the disk $\Delta$ under~$f$}.
\end{defn}

\begin{remarks}
\mbox{}
\begin{enumerate}[a)]
\item Note that condition~b) implies that $\alpha([0,1])$ intersects
$P\cup Q$ only at $p$ and~$q$.
\item Condition~a) is ponderous, but the reason for it is simply
  explained. In the applications in this paper, where~$f$ is the
  horseshoe, the arc~$\alpha$ will be taken to be an arc of the stable
  manifold of~$p$ followed by an arc of the unstable manifold
  of~$q$. Suppose, for example, that $p$~is a periodic point of
  negative index (so that $f^n$ sends each branch of its stable
  manifold to itself), while $q$~is a periodic point of positive index
  (so that~$f^n$ sends each branch of its unstable manifold to the
  other branch). The parameter~$b$ corresponds to the point~$x$ where
  the two manifolds meet, so that $f^n(\alpha([0,b]))$ is a subset
  $\alpha([0,a])$ of the image of~$\alpha$: the remainders
  $\alpha([a,1])$ and $f^n(\alpha([b,1]))$ of the arc and its image
  under~$f^n$ are disjoint except at their endpoints (which are at~$x$
  and~$q$), and hence form a simple closed curve. See
  Figure~\ref{fig:gammaconfig}, for example.
\end{enumerate}
\end{remarks}

\begin{thm}
\label{thm:arc}
Let $f\colon D^2\to D^2$ be an orientation-preserving homeomorphism,
$p$ and $q$ be distinct period $n$ points of $f$, and $R$ be a finite
subset of $D^2$ with $f(R)=R$. Suppose that there exists an arc
$\alpha\colon[0,1]\to D^2\setminus R$ from $p$ to $q$, which defines a
disk~$\Delta$ under~$f$ which is disjoint from~$R$. Then $(p\,;\,f)$
and $(q\,;\,f)$ are connected by isotopy rel. $R$.
\end{thm}

\begin{proof}
Denote by $\alpha_i$ the image of the arc $f^i\circ\alpha$ for $0\le
i\le n$: since $f$ is a homeomorphism, it follows from
Definition~\ref{defn:define}~b) that $\alpha_i\cap\alpha_j=\emptyset$
for $i\not=j$, except if $\{i,j\}=\{0,n\}$.

Thus $\alpha_i$ has endpoints on $P\cup Q$ and is disjoint from
$\partial\Delta$ for $1\le i<n$, hence by
Definition~\ref{defn:define}~c) $\alpha_i\cap\Delta=\emptyset$ for
$1\le i<n$. Let $D$ be a disk which contains $\alpha_0\cup\alpha_n$ in
its interior, but which is disjoint from $R$ and from $\alpha_i$ for
$1\le i<n$; and let $C$ be a simple closed curve bounding a disk
containing $\alpha_0$ in its interior such that both $C$ and $f^n(C)$
are contained in $\Int\,D$, while $f^i(C)$ is disjoint from $D$ for
$1\le i<n$.

Then $f^n(C)$ is isotopic to $C$ rel. $P\cup Q\cup R$, since both are
simple closed curves contained in~$D$, surrounding the only two points
of $P\cup Q\cup R$ which lie in~$D$. Hence (by a theorem of
Epstein~\cite{Ep}) there is a homeomorphism $h\colon D^2\to D^2$,
supported in $D$ and isotopic to the identity rel. $P\cup Q\cup R$,
with $h(f^n(C))=C$. Let $F=h\circ f$, so that $F^n(C)=C$ and $f\simeq
F$ rel. $P\cup Q\cup R$. Clearly $(p\,;\,f)\sim(p\,;\,F)$ and
$(q\,;\,f)\sim(q\,;\,F)$ (using the constant paths from $p$ to $p$ and
from $q$ to $q$, and the isotopy $f\simeq F$ rel. $P\cup Q\cup R$), so
it remains to show that $(p\,;\,F)\sim(q\,;\,F)$.

Let $E_0$ be the closed disk bounded by $C$, write $E_i=F^i(E_0)$ for
$1\le i<n$, and set $E=\bigcup_{i=0}^{n-1}E_i$. Thus $E$ is an
$F$-invariant subset of $D^2$, disjoint from $R$, consisting of $n$
mutually disjoint disks. Let $H\colon D^2\to D^2$ be a homeomorphism
supported in $E$ with $H(F^i(q))=F^i(p)$ for each $i$, and let
$\{H_t\}\colon H\simeq \id$ be an isotopy supported in~$E$. Then
$H_t(q)$ is a period $n$ point of $H_t\circ F\circ H_t\I$ for each
$t$, so that $(p\,;\,G)\sim(q\,;\,F)$, where $G=H\circ F\circ
H\I$. However $G$ agrees with $F$ on $P$ and outside of $E$: thus
applying the Alexander trick to all of the components of $E$ (which
each contain a single point of $P$), there is an isotopy from $G$ to
$F$ relative to $P\cup R$, which provides an isotopy connection
$(p\,;\,F)\sim (p\,;\,G)$. Hence $(p\,;\,F)\sim(q\,;\,F)$ as required.
\end{proof}

\subsection{Pruning theory}
\label{sec:pruning}
Pruning theory provides a means of destroying some of the dynamics of
a surface homeomorphism by an isotopy with controlled support. The
following definition and theorem are from~\cite{Pruning}, simplified
in accordance with the requirements of this paper.

\begin{defns}
\label{defn:pruning}
A {\em pruning disk} for the horseshoe map $F:D^2\to D^2$ is a closed
topological disk $\Delta\subseteq D^2$ whose boundary is the union of an
arc~$C$ of stable manifold and an arc~$E$ of unstable manifold,
intersecting only at their endpoints, which satisfy
\[F^n(C)\cap\Int(\Delta) = F^{-n}(E)\cap\Int(\Delta)=\emptyset \qquad\text{for
  all }n\ge 1.\]
The common endpoints of the two arcs are called the {\em vertices} of
$\Delta$. 
\end{defns}

\begin{thm}
\label{thm:pruning}
Let~$\Delta$ be a pruning disk for~$F$. Then there exists an isotopy,
supported in $\bigcup_{n\in\Z}F^n(\Int(\Delta))$, from~$F$ to a
homeomorphism~$F_\Delta$ for which all points of $\Int(\Delta)$ are wandering.
\end{thm}

The pruning isotopy therefore destroys all of the dynamics in
$\Int(\Delta)$, while leaving untouched any orbits which do not enter
$\Int(\Delta)$. In this paper the pruning theorem will be applied in the
form of the following corollary:
\begin{cor}
\label{cor:pruning}
Let~$w$ be a lone decoration and $q\in(0,q_w)$. Suppose that there is
a pruning disk~$\Delta$ containing points of all four periodic orbits
$P^w_q$ in its interior. If~$R$ is a horseshoe periodic orbit disjoint
from~$\Delta$, then $R\not\ge P^w_q$.
\end{cor}
\begin{proof}
The set of braid types of the pruned homeomorphism~$F_\Delta$ is precisely
the set of braid types of periodic orbits of~$F$ which are disjoint
from~$\Delta$. In particular, $F_\Delta$ has a periodic orbit of the braid type
of~$R$, but none of the braid type of~$P^w_q$.
\end{proof}

\section{Proof of the main theorem}
\label{sec:proof}
Let $w$ be a lone decoration of length~$k$. This decoration will be
fixed throughout the section, and therefore the dependence of many
objects on~$w$ will be suppressed: on the few occasions when it is
temporarily important to indicate this dependence, this will be done
by means of a superfix~$w$. For the sake of clarity, it is assumed at
first that $w$ is even: the modifications necessary in the case of
odd~$w$ are described in Section~\ref{sec:odd}.

\subsection{Iterated arcs}
\label{sec:itarc}
The proof of the theorem depends on the details of the
configuration of the $F$-images of a collection of arcs joining points
of the four periodic orbits of height $q$ and decoration $w$ for each
$q\in(0,q_w)\cap\Q$. For the remainder of this subsection, let $q=m/n$
be a fixed element of $(0,q_w)$. For the sake of notational clarity,
arcs $\alpha:[0,1]\to D^2$ and their images $\alpha([0,1])$ will not
be distinguished carefully; and points $x\in\Omega(F)$ will often be
identified with $k(x)\in\{0,1\}^\Z$.

Let $p_1=\ol{c_q0w0}$, $p_2=\ol{c_q0w1}$, $p_3=\ol{c_q1w1}$, and
$p_4=\ol{c_q1w0}$ be the rightmost points on each of the four periodic
orbits of height $q$ and decoration $w$ (which have
period~$n+k+3$). Note that $\Index(p_i,F^{n+k+3})=(-1)^i$, since the words
$c_q0w0$ and $c_q1w1$ are even and the words $c_q0w1$ and $c_q1w0$
are odd.

\begin{figure}[htbp]
\lab{1}{p_1}{t}
\lab{2}{p_4}{t}
\lab{9}{p_2}{l}
\lab{8}{p_3}{l}
\lab{7}{\delta}{l}
\lab{3}{\gamma}{r}
\lab{4}{S}{l}
\lab{C}{C}{l}
\begin{center}
\pichere{0.4}{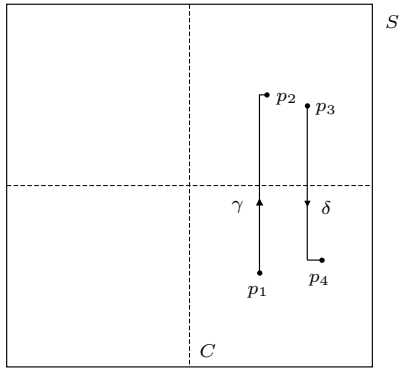}
\end{center}
\caption{The arcs $\gamma$ and $\delta$}
\label{fig:arcs}
\end{figure}

Let $\gamma$ and $\delta$ be the following arcs connecting these
points:

\begin{enumerate}[a)]
\item $\gamma$ goes from $p_1$ to $p_2$. It is the concatenation of
  the vertical arc from $p_1$ to the point
  $\lib{c_q0w1}\cdot\rib{c_q0w0}$ and the horizontal arc from
  $\lib{c_q0w1}\cdot\rib{c_q0w0}$ to~$p_2$.
\item $\delta$ goes from $p_3$ to $p_4$. It is the concatenation of
  the vertical arc from $p_3$ to the point
  $\lib{c_q1w0}\cdot\rib{c_q1w1}$ and the horizontal arc from
  $\lib{c_q1w0}\cdot\rib{c_q1w1}$ to~$p_4$.
\end{enumerate}
These arcs are depicted schematically in
Figure~\ref{fig:arcs}. The
points~$p_i$ are shown with their correct relative horizontal and
vertical orderings, calculated using the fact that $w$ is even. 

It will be shown that each of these arcs defines a disk (in the sense
of Definition~\ref{defn:define}) under the horseshoe map~$F$, and
explicit conditions for a point of the non-wandering set $\Omega(F)$
to lie in each of the disks thus defined will be given.

\begin{figure}
\lab{1}{p_1}{t}
\lab{9}{p_2}{b}
\lab{3}{\gamma}{r}
\lab{4}{S}{l}
\lab{2}{\gamma_{n+k+3}}{t}
\lab{C}{C}{l}
\begin{center}
\pichere{0.4}{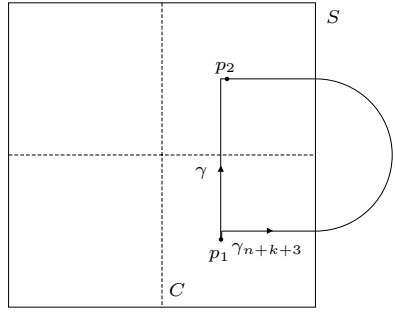}
\end{center}
\caption{The disk defined by~$\gamma$}
\label{fig:gammaconfig}
\end{figure}

\begin{lem}
\label{lem:gamma}
The arc~$\gamma$ defines a disk~$C_q$ under~$F$. A point $x=b\cdot f$ of
$\Omega(F)$ lies in $\Int\,C_q$ if and only if
\begin{enumerate}[a)]
\item $f\succ \rib{c_q0w0}$, and
\item $\sigma(b)\succ \rib{\hw 0c_q1}$.
\end{enumerate}
\end{lem}
\begin{proof}
Denote by $\gamma_j$ the arc $F^j\circ\gamma$, for $0\le j\le
n+k+3$. A straightforward induction shows that for $0\le j\le n+k+2$,
the arc $\gamma_j$ is contained in~$S$, and is the concatenation of
the vertical arc from $\sigma^j(\ol{c_q0w0})$ to
$\sigma^j(\lib{c_q0w1}\cdot\rib{c_q0w0})$ and the horizontal arc from
$\sigma^j(\lib{c_q0w1}\cdot\rib{c_q0w0})$ to
$\sigma^j(\ol{c_q0w1})$. This is true by definition when~$j=0$, and
follows for $1\le j\le n+k+2$ since $\gamma_{j-1}$ does not cross the
vertical centre line~$C$ of~$S$.

Thus $\gamma_{n+k+2}$ is the concatenation of the vertical arc from
$\ol{0c_q0w}$ to $\lib{1c_q0w}\cdot\rib{0c_q0w}$ and the horizontal
arc from $\lib{1c_q0w}\cdot\rib{0c_q0w}$ to $\ol{1c_q0w}$ (which
crosses~$C$). Hence (see Figure~\ref{fig:gammaconfig})
$\gamma_{n+k+3}$ is the concatenation of
\begin{enumerate}[i)]
\item the vertical arc from $p_1=\ol{c_q0w0}$ to
$\lib{1c_q0w}0\cdot\rib{c_q0w0}$;
\item the horizontal arc from $\lib{1c_q0w}0\cdot\rib{c_q0w0}$ to
the right hand edge of~$S$;
\item a semicircular arc outside of~$S$; and
\item the horizontal arc from the right hand edge of~$S$ to
$\ol{c_q0w1}=p_2$.
\end{enumerate}
Condition~a) of Definition~\ref{defn:define} thus follows (with
$\alpha(b)=\lib{1c_q0w}0\cdot\rib{c_q0w0}$) and conditions~b) and~c)
are satisfied because~$p_1$ and~$p_2$ lie to the right of all other
points in both of their orbits, and the arcs $\gamma_j$ are contained
in~$S$ for $1\le j\le n+k+2$. A point $x=b\cdot f$ of $\Omega(F)$ lies
in $\Int\,C_q$ if and only if both $f\succ\rib{c_q0w0}$, and
\mbox{$0\rib{\hw 0c_q1}\prec b\prec \rib{1\hw 0c_q}$}: the latter
condition is equivalent to $\sigma(b)\succ \rib{\hw 0c_q1}$.
\end{proof}

The corresponding result for the arc~$\delta$ is analogous, and its
proof is omitted.

\begin{lem}
\label{lem:delta}
The arc~$\delta$ defines a disk~$D_q$ under~$F$. A point $x=b\cdot f$ of
$\Omega(F)$ lies in $\Int\,D_q$ if and only if
\begin{enumerate}[a)]
\item $f\succ \rib{c_q1w1}$, and
\item $\sigma(b)\succ \rib{\hw 1c_q0}$.
\end{enumerate}
\end{lem}

For the purposes of the current argument, denote the arc~$\gamma$,
which depends on~$q$ and on the decoration~$w$, by $\gamma^w_q$. Now
$q<q_w=q_{\hat{w}}$ (Lemma~\ref{lem:revq}), and hence, applying
Lemma~\ref{lem:gamma} with decoration $\hat{w}$, the arc
$\gamma^{\hat{w}}_q$ (which joins $\overline{c_q0\hat{w}0}$ to
$\overline{c_q0\hat{w}1}$) defines a disk $C_q^{\hat{w}}$ under~$F$,
and a point $x=b\cdot f$ of $\Omega(F)$ lies in $\Int C_q^{\hat{w}}$
if and only if
\begin{enumerate}[a)]
\item $f\succ (c_q0\hat{w}0)^\infty$, and
\item $\sigma(b)\succ(w0c_q1)^\infty$.
\end{enumerate}

Applying the involution $\phi$ of Section~\ref{sec:hsinverse}, the arc
$\phi(\gamma^{\hat{w}}_q)$ (which joins $\overline{0w0c_q}$ to
$\overline{1w0c_q}$) defines a disk $\phi(C_q^{\hat{w}})$
under~$F^{-1}$, and a point $x=b\cdot f$ of $\Omega(F)$ lies in
$\Int\phi(C_q^{\hat{w}})$ if and only if
\begin{enumerate}[a)]
\item $b\succ (c_q0\hat{w}0)^\infty$, and
\item $\sigma(f)\succ(w0c_q1)^\infty$.
\end{enumerate}

\begin{figure}[htbp]
\lab{1}{r_1}{t}
\lab{2}{r_2}{b}
\lab{9}{r_4}{b}
\lab{8}{r_3}{l}
\lab{7}{\beta}{b}
\lab{3}{\alpha}{t}
\lab{4}{S}{l}
\lab{C}{C}{l}
\begin{center}
\pichere{0.4}{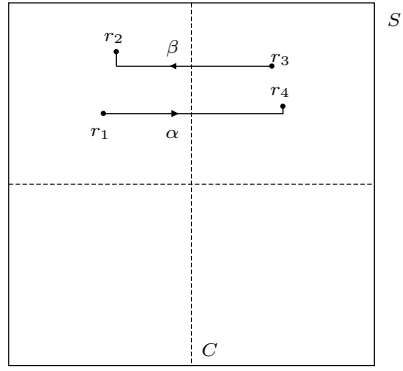}
\end{center}
\caption{The arcs $\alpha$ and $\beta$}
\label{fig:arcs2}
\end{figure}

Write $\alpha=\alpha^w_q=\phi(\gamma^{\hat{w}}_q)$, and similarly
$\beta=\beta^w_q=\phi(\delta^{\hat{w}}_q)$ (see
Figure~\ref{fig:arcs2}). Thus $\alpha$ joins the highest
point~$r_1=\overline{0w0c_q}$ on the orbit of~$p_1$ to the highest
point $r_4=\overline{1w0c_q}$ on the orbit of~$p_4$; while $\beta$
joins the highest point~$r_2=\overline{0w1c_q}$ on the orbit of~$p_2$
to the highest point~$r_3=\overline{1w1c_q}$ on the orbit of~$p_3$,
and the argument above gives:
\begin{lem}
\label{lem:alpha}
The arc~$\alpha$ defines a disk~$A_q$ under~$F^{-1}$. A
point~$x=b\cdot f$ of $\Omega(F)$ lies in $\Int A_q$ if and only if
\begin{enumerate}[a)]
\item $b\succ (c_q0\hat{w}0)^\infty$, and
\item $\sigma(f)\succ(w0c_q1)^\infty$.
\end{enumerate}
\end{lem}

\begin{lem}
\label{lem:beta}
The arc~$\beta$ defines a disk~$B_q$ under~$F^{-1}$. A
point~$x=b\cdot f$ of $\Omega(F)$ lies in $\Int B_q$ if and only if
\begin{enumerate}[a)]
\item $b\succ (c_q1\hat{w}1)^\infty$, and
\item $\sigma(f)\succ(w1c_q0)^\infty$.
\end{enumerate}
\end{lem}

\begin{figure}[htbp]
\begin{center}
\pichere{0.5}{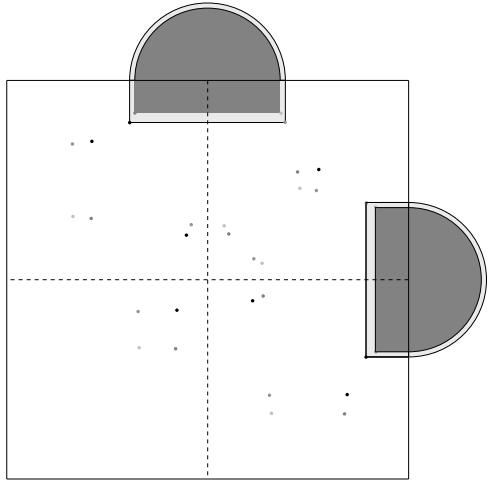}
\end{center}
\caption{The disks $A_{1/3}$, $B_{1/3}$, $C_{1/3}$, and $D_{1/3}$ for $w=11$}
\label{fig:alldisk}
\end{figure}

Figure~\ref{fig:alldisk} depicts the disks $A_{1/3}$, $B_{1/3}$,
$C_{1/3}$, and $D_{1/3}$ for $w=11$, together with the four periodic
orbits of height~$1/3$ and decoration~$w$. This figure is drawn to
scale, and short segments of stable and unstable manifolds cannot be
discerned on it. The disks $A_{1/3}$ and $C_{1/3}$ are shaded lightly,
and the disks $B_{1/3}$ and $D_{1/3}$ are shaded heavily.

This figure motivates the following straightforward lemma.

\begin{lem}
\label{lem:contained}
For any decoration~$w$ and any~$q\in(0,q_w)$, $D_q\subseteq
C_q$ and $B_q\subseteq A_q$.
\end{lem}

\begin{proof}
Let $x=b\cdot f\in\Omega(F)$. It suffices to show that if
$x\in\Int(D_q)$ then \mbox{$x\in\Int(C_q)$}. That $B_q\subseteq A_q$
will then follow since $A_q=\phi(C_q^{\hat{w}})$ and
$B_q=\phi(D_q^{\hat{w}})$.

Suppose $x\in \Int(D_q)$: then $f\succ \rib{c_q1w1}$ and
$\sigma(b)\succ \rib{\hat{w}1c_q0}$ by Lemma~\ref{lem:delta}. Since
both $c_q$ and $\hat{w}$ are even words, it follows that
$f\succ\rib{c_q0w0}$ and \mbox{$\sigma(b)\succ\rib{ \hat{w}0c_q1}$},
so that $x\in \Int(C_q)$ by Lemma~\ref{lem:gamma}.
\end{proof}

\subsection{Linking properties}
Let $q\in(0,q_w)\cap\Q$, and for $1\le i\le 4$ denote by $P_{i,q}$ the
periodic orbit containing the point~$p_i$ of
Section~\ref{sec:itarc}. Fix a suspension $\{f_t\}$ of the horseshoe
map~$F$. Given a horseshoe periodic orbit~$R$ distinct from
each~$P_{i,q}$, denote by~$L_{i,q}(R)$ the linking number of~$P_{i,q}$
about~$R$ with respect to the suspension~$\{f_t\}$. 

\begin{lem}
\label{lem:linking}
\begin{eqnarray*}
L_{4,q}(R)&=&L_{1,q}(R)+|R\cap A_q|\\
L_{3,q}(R)&=&L_{2,q}(R)+|R\cap B_q|\\
L_{2,q}(R)&=&L_{1,q}(R)+|R\cap C_q|\\
L_{3,q}(R)&=&L_{4,q}(R)+|R\cap D_q|.
\end{eqnarray*}
In particular, 
\begin{eqnarray*}
\max(L_{1,q}(R),L_{2,q}(R),L_{3,q}(R),L_{4,q}(R))&=&L_{3,q}(R)
\qquad\text{ and}\\
\min(L_{1,q}(R),L_{2,q}(R),L_{3,q}(R),L_{4,q}(R))&=&L_{1,q}(R).
\end{eqnarray*}
\end{lem}
\begin{proof}
Note first that the boundaries of the disks $A_q$, $B_q$, $C_q$, and
$D_q$ are composed of segments of the stable and unstable manifolds of
points of the orbits $P_{i,q}$, and so cannot intersect the orbit~$R$.

Consider the arc $\gamma$ from $p_1$ to $p_2$. The curve
$\gamma\cdot(F^{n+k+3}(\gamma))^{-1}$ goes clockwise around the
boundary of~$C_q$. Hence $L_{1,q}(R)=L_{2,q}(R)-|R\cap C_q|$ by
Lemma~\ref{lem:linkcalc} as required. Similarly
$L_{3,q}(R)=L_{4,q}(R)+|R\cap D_q|$ (the curve
$\delta\cdot(F^{n+k+3}(\delta))^{-1}$ goes anti-clockwise around the
boundary of~$D_q$).

Let $\tilde\alpha_q = F^{-(n+k+3)}(\alpha_q)$, an arc from $r_1$ to
$r_4$. Then $\tilde\alpha_q\cdot \alpha_q^{-1}$ goes clockwise around
the boundary of~$A_q$, so that $L_{1,q}(R)=L_{4,q}(R)-|R\cap A_q|$ as
required. Similarly $L_{3,q}(R)=L_{2,q}(R)+|R\cap B_q|$.
\end{proof}

\begin{cor}
\label{cor:linkcoincidence}
Let~$R$ be any horseshoe periodic orbit distinct from
each~$P_{i,q}$. Then
\[|R\cap C_q|+|R\cap B_q|=|R\cap D_q|+|R\cap A_q|.\]
\end{cor}

\begin{proof}
Both are equal to $L_{3,q}(R)-L_{1,q}(R)$.
\end{proof}

\subsection{Forcing conditions}
The first main result is that $\{P_q^w\,:\,q\in\Q^w\}$ is totally
ordered by the forcing relation. The following lemma will be used:

\begin{lem}
\label{lem:Cq}
Let $q,q'\in (0,q_w)\cap\Q$ with $q'<q$. Let $R$ be the periodic orbit
of code $c_{q'}0w0$. Then $|R\cap A_q|=|R\cap C_q|=1$, and $|R\cap
B_q|=|R\cap D_q|=0$.
\end{lem}

\begin{proof}
Observe first that a point of~$R$ lies in one of the disks if and only
if it lies in its interior, since~$R$ is distinct from the
orbits~$P_{i,q}$. 

\begin{enumerate}[a)]
\item The rightmost point~$x=b\cdot f$ of~$R$ lies in $C_q\setminus
  D_q$. For $f=\rib{c_{q'}0w0}$ and $b=\rib{0\hat{w}0c_{q'}}$. Now
  $q'=q(f)<q=q(\rib{c_q0w0})$ and hence $f\succ\rib{c_q0w0}$, so
  condition~a) of Lemma~\ref{lem:gamma} (for a point to lie in~$C_q$)
  is satisfied.

Since $q_w=q_{\widehat{w}}$ (Lemma~\ref{lem:revq}) and $q'<q<q_w$, the
words $c_{q'}0\hat{w}0$ and $c_q1\hat{w}0$ are codes of periodic
orbits of heights $q'$ and $q$ by Lemma~\ref{lem:qw}. Therefore
\mbox{$\rib{c_{q'}0\hat{w}0}\succ\rib{c_q1\hat{w}0}$}. Prepending the
even word $\hat{w}0$ to both sides gives
\mbox{$\sigma(b)\succ\rib{\hat{w}0c_q1}$}, so that condition~b) of
Lemma~\ref{lem:gamma} is satisfied.

On the other hand, $\sigma(b)=\hat{w}0\ldots \prec\rib{\hat{w}1c_q0}$
(since $\hat{w}$ is even), and hence condition~b) of
Lemma~\ref{lem:delta} (for a point to lie in~$D_q$) is not satisfied. 
\item No other point $x=b\cdot f$ of~$R$ lies in~$C_q$.  Certainly
  this is true if~$f$ begins in~$w$ or at the final~$1$ of~$c_q'$: for
  \[q_w=\min_{0\le i\le k+2}q\left(\sigma^i(\rib{10w0})\right)\]
gives that $q(f)\ge q_w>q=q(\rib{c_q0w0})$ for any $f\in\{0,1\}^\N$
which starts either with $10w010$, or with a final subword of~$w$
followed by $010$, so that $f\prec\rib{c_q0w0}$, and condition~a) of
Lemma~\ref{lem:gamma} is not satisfied. Thus condition~a) of this
lemma can only be satisfied if
\begin{eqnarray*}
f&=&10^{\ka_i}1^2\ldots 1^20^{\ka_{m'}}10w0\ldots\qquad\text{and}\\
b&=&10^{\ka_{i-1}}1^20^{\ka_{i-2}}1^2\ldots
1^20^{\ka_1}10\hat{w}0\ldots \qquad\text{for some~$i$ with $1<i\le m'$}
\end{eqnarray*}
(where $c_q' = 10^{\kappa_1}1^2\ldots 1^20^{\kappa_{m'}}1$). 
Hence
\[q(10\sigma(b))\le q'<q_w=q_{\hat{w}}\le q\left(10\rib{\hat{w}0c_q1}
\right)\] (the first inequality coming from Lemma~\ref{lem:starlem}
and the fact that $c_{q'}$ is palindromic, and the second from the
above formula for $q_w$), so that $10\sigma(b)\succ
10\rib{\hat{w}0c_q1}$, i.e. \mbox{$\sigma(b)\prec \rib{\hat{w}0c_q1}$}, and
condition~b) of Lemma~\ref{lem:gamma} fails.
\end{enumerate}
Hence $|R\cap C_q|=1$ and $|R\cap D_q|=0$. So $|\widehat{R}\cap
C_q^{\hat{w}}|=1$ and $|\widehat{R}\cap D_q^{\hat{w}}|=0$, and
therefore $|R\cap A_q|=|\phi(\widehat{R})\cap\phi(C_q^{\hat{w}})|=1$
and $|R\cap B_q|=|\phi(\widehat{R})\cap\phi(D_q^{\hat{w}})|=0$ as
required. 
\end{proof}

Theorem~\ref{thm:main}~a) follows immediately from the following
result.

\begin{thm}
\label{thm:order}
Let $q\in\Q^w$ and $q'\in(0,q)\cap\Q$. Then $P_{q'}^w\ge P_q^w$.
\end{thm}
\begin{proof}
Let~$R$ be the periodic orbit of code $c_{q'}0w0$.
Lemmas~\ref{lem:linking} and~\ref{lem:Cq} give
that \[L_{2,q}(R)=L_{3,q}(R)=L_{4,q}(R)=L_{1,q}(R)+1.\]
 
Suppose first that $L_{1,q}(R)$ is not divisible by $n+k+3$. It will
be shown that $(p_1\,;\,F)$ is unremovable in~$(D^2,R)$, which will
establish the result. Note first that $(p_1\,;\,F)$ is certainly
separated from~$R$, since~$R$ and~$P_{1,q}$ have different braid types
(Remark~\ref{rem:separated}). To show that $(p_1\,;\,F)$ is
uncollapsible, suppose that $g_j\to g$ is a sequence of homeomorphisms
$D^2\to D^2$, and $r_j\to r$ in $D^2$ is such that $r_j$ is a
period~$n+k+3$ point of $g_j$ with $(r_j\,;\,g_j)\sim(p_1\,;\,F)$ for
each~$j$. The fact that $n+k+3\not|L_{1,q}$ means (by
Lemma~\ref{lem:nocollapse}) that~$r$ cannot be a fixed point
of~$g$. If~$r$ were a period $(n+k+3)/\ell$ point of~$g$ for some
$\ell$ with $1<\ell<n+k+3$, then for $j$ sufficiently large the braid
type of the $g_j$-orbit of~$r_j$ would be reducible
with~$(n+k+3)/\ell$ reducing curves, contradicting the fact that
$P_{1,q}$ has pseudo-Anosov braid type. Thus~$r$ is a period~$n+k+3$
point of~$g$, establishing the uncollapsibility of
$(p_1\,;\,F)$. Finally, to show that $(p_1\,;\,F)$ is essential: if
$r\in\snc(p_1\,;\,F)$, then~$r$ lies on a periodic orbit of braid type
$\bt(P_{1,q})$ and has linking number~$L_{1,q}(R)$ about~$R$ with
respect to some suspension of~$F$ by
Lemma~\ref{lem:icpreserve}. Since~$w$ is a lone decoration and
$L_{2,q}(R)=L_{3,q}(R)=L_{4,q}(R)=L_{1,q}(R)+1$, it follows that $r\in
P_{1,q}$, and hence $I(p_1\,;\,F)=-|\snc(p_1\,;\,F)|\not=0$, so
$(p_1\,;\,F)$ is essential as required.

 On the other hand, if $L_{1,q}(R)$ is divisible by $n+k+3$, then
 $L_{2,q}(R)$, $L_{3,q}(R)$, and $L_{4,q}(R)$ are not. The above
 argument can then be repeated with~$p_2$, the only part which is any
 different being the proof that $I(p_2\,;\,F)\not=0$. For this part,
 observe first that by Lemma~\ref{lem:icpreserve} the only points
 which can lie in $\snc(p_2\,;\,F)$ are points of the form $F^i(p_j)$
 for $j\in\{2,3,4\}$, since these are the only periodic points of~$F$
 which lie on periodic orbits of braid type $\bt(P_{2,q})$ and have
 linking number $L_{2,q}(R)$ about~$R$. However, since $|R\cap
 B_q|=|R\cap D_q|=0$, Theorem~\ref{thm:arc} gives that
 $(p_2\,;\,F)\sim(p_3\,;\,F)\sim(p_4\,;\,F)$. Hence if
 $F^i(p_j)\in\snc(p_2\,;\,F)$ for some $j\in\{2,3,4\}$, then
 (Lemma~\ref{lem:iciterate}) this is true for all
 $j\in\{2,3,4\}$. Since points of $P_{2,q}$ and $P_{4,q}$ have index
 $+1$, while those of $P_{3,q}$ have index $-1$, it follows that
 $I(p_2\,;\,F)=|\snc(p_2\,;\,F)|/3 > 0$ as required.
\end{proof}

\begin{remark}
  A cleaner approach to this proof would be to start by using the fact
  that $R\cap D_q = \emptyset$ to prune away the periodic orbits
  $P_{3,q}$ and $P_{4,q}$ by an isotopy which leaves~$R$, $P_{1,q}$,
  and $P_{2,q}$ untouched (a suitable pruning disk can be constructed
  by analogy with Lemma~\ref{lem:enlarge} below). There would then remain
  only two periodic orbits of the braid type of $P_{1,q}$, whose
  linking numbers about~$R$ would differ by~$1$. However, in order to
  take this approach it is necessary to ensure that the indices of
  $p_1$ and $p_2$ are unchanged after the pruning isotopy, which
  requires more careful control of the support of this isotopy. The
  technical details involved in doing this are more complicated than
  the approach taken in the proof above.
\end{remark}

\begin{figure}[htbp]
\begin{center}
\lab{1}{p_1}{t}
\lab{2}{p_4}{t}
\lab{9}{p_2}{l}
\lab{8}{p_3}{l}
\lab{7}{\delta}{l}
\lab{3}{\gamma}{r}
\lab{4}{S}{l}
\lab{A}{v_1}{r}
\lab{B}{v_0}{r}
\lab{C}{\Delta_q}{bl}
\lab{E}{E}{tl}
\lab{F}{C}{r}
\pichere{0.6}{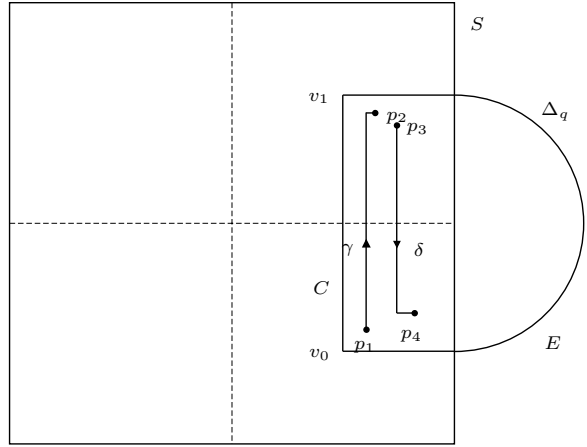}
\end{center}
\caption{The pruning disk $\Delta_q$}
\label{fig:delta}
\end{figure}

The next lemma describes how the disk $C_q$ can be enlarged to a
pruning disk $\Delta_q$ which contains all of the points $p_i$. (See
Figure~\ref{fig:delta}.) $\Delta_q$ is ``approximately the same'' as
$C_q$ in the sense that if $q$ has large denominator then
the boundary of $\Delta_q$ is close to that of~$C_q$, and in
particular $\Delta_q\setminus C_q$ cannot contain periodic points of
low period. This idea is encapsulated in the following definition:
\begin{defn}
\label{defn:approx}
Let $q=m/n\in\Q$. If $X,Y\subseteq D^2$, write $X\subseteq_q Y$ if any
periodic point of~$F$ of period less than $n/2$ which lies in~$X$ also
lies in~$Y$.
\end{defn}

\begin{lem}
\label{lem:enlarge}
Let~$v_0$ and $v_1$ be the points
$\li{}\,0c_q0w\opti\cdot\rib{c_q0w011}$. Then $v_0$ and $v_1$ are the
vertices of a pruning disk $\Delta_q$ which
contains~$\{p_1,p_2,p_3,p_4\}$. Moreover $\Delta_q\subseteq_q C_q$.
\end{lem}
\begin{proof}
Observe that a point $x=b\cdot f$ of $\Omega(F)$ lies in
$\Int(\Delta_q)$ if and only if \mbox{$f\succ\rib{c_q0w011}$} and
$\sigma(b)\succ\hat{w}0c_q\ri{0}$.

Let~$C$ and~$E$ be the segments of stable and unstable manifold
constituting the boundary of~$\Delta_q$. Observe that, for $n\ge 1$,
$F^n(C)$ (respectively $F^{-n}(E)$) is a segment of stable
(respectively unstable) manifold, with endpoints $F^n(v_i)$
(respectively $F^{-n}(v_i)$), which is contained in the central
square~$S$ of the horseshoe. Thus to show that $\Delta_q$ is a pruning
disk (Definition~\ref{defn:pruning}), it is enough to show that
$F^n(v_i)\not\in\Int(\Delta_q)$ for all $n\in\Z$.

Now $q_{w01}=\min_{0\le i\le k+4}q(\sigma^i(\rib{10w010})) =\min_{0\le
  i\le k+2}q(\sigma^i(\rib{10w0}))=q_w$, and hence $q<q_{w01}$. Thus
(Lemma~\ref{lem:qw}) $c_q0w011$ is the code of a periodic orbit, whose
rightmost point is on the stable leaf containing~$C$. It follows that
$F^n(v_i)\not\in\Int(\Delta_q)$ for all~$n\ge 0$.

If $n<0$, then $h(\sigma^n(v_i))\succ\rib{c_q0w011}$ could only occur
if 
\[h(\sigma^n(v_i)) = 10^{\kappa_r}1^2\ldots
1^20^{\kappa_m}10w\opti\rib{c_q0w011}\] for some $1\le r\le m$ (recall
the horiztonal and vertical coordinate functions~$h$ and~$v$ from
Section~\ref{sec:smal-hors-unim}). If $r>1$ then this is not possible
(as $q(10^{\kappa_r}1^2\ldots 1^20^{\kappa_m}10\ldots) > q$), while if
$r=1$ then $\sigma(v(\sigma^n(v_i)))=0^\infty\not\succ
\hat{w}0c_q0^\infty$. Thus $\Delta_q$ is a pruning disk as required.

It is straightforward that $\Delta_q$ contains all of the
$p_i$, since \mbox{$\rib{c_q\opti w\opti}\succ\rib{c_q0w011}$} and
$\rib{\hat{w}\opti c_q\opti} \succ \hat{w}0c_q\ri{0}$.

For the final part of the lemma, suppose that $x=\overline{c}$ is a
point of a period \mbox{$N=|c|<n/2$} orbit~$R$ of~$F$. It is required
to show that if $x\in\Int(\Delta_q)$ then
\mbox{$x\in\Int(C_q)$}. Suppose for a contradiction, then, that
$x\in\Int(\Delta_q)\setminus\Int(C_q)$: that is, by
Lemma~\ref{lem:gamma}, either
\[
\begin{array}{rcccl}
\rib{c_q0w011} & \prec & \ri{c} & \preceq &
\rib{c_q0w0},\quad\text{or}\\
\hat{w}0c_q\ri{0} & \prec & \sigma(\ri{\hat{c}}) & \preceq &
\rib{\hat{w}0c_q1}.
\end{array}
\]
If the former inequalities hold, then $q(c^\infty)=q$. Since
$|c|<n/2$, both $c$ and $cc$ must be initial words of $c_q$, and so
$c=10^{\kappa_1}1^2\ldots 1^20^{\kappa_r}1$ for some $r < m/2$. Now
the code~$c_R$ of~$R$ is some cyclic permutation of $c$, so $c_R =
10^{\kappa_i}1^2 \ldots 1^20^{\kappa_r}1^20^{\kappa_1}\ldots 1^2
0^{\kappa_{i-1}}1$ for some $i$ between $1$ and~$r$. Removing the even
initial subword $10^{\kappa_1}1^2\ldots 1^20^{\kappa_{i-1}}1$ (of
length~$s$, say) from each term in the inequalities gives
\[\sigma^s (\rib{c_q0w011}) \prec \ri{c_R} \preceq
\sigma^s(\rib{c_q0w0}),\]
and hence $q(R)=q(\ri{c_R})\ge q$. On the other hand, $q(R)\le
q(c^\infty)=q$, and hence $q(R)=q=m/n$, and so~$R$ has period at
least~$n$. This is the required contradiction.

If the latter inequalities hold, then (removing the even initial
subword $\hat{w}0$ from each term) there is some point $d^\infty$
of~$R$ such that
\[c_q0^\infty \prec \ri{\hat{d}} \preceq \rib{c_q1\hat{w}0},\]
and the proof proceeds in exactly the same way.
\end{proof}

The description of the invariant $r^w$ follows relatively easily from
the following theorem:

\begin{thm}
\label{thm:forcing}
Let~$R$ be a period~$N$ horseshoe orbit, and let
$q=m/n\in\Q^w$ be such that $n>2N$. Then $R\ge P_q^w$ if and only if
$R\cap A_q\not=\emptyset$ and $R\cap C_q\not=\emptyset$.
\end{thm}

\begin{proof}
Notice that~$R$ cannot intersect the boundary of any of the
disks~$A_q$, $B_q$, $C_q$, $D_q$, since its period is less than the
period $n+k+3$ of $P_q^w$.

If $R\cap C_q = \emptyset$, then $R\not\ge P^w_q$ by
Lemma~\ref{lem:enlarge} and Corollary~\ref{cor:pruning}. Similarly, if
$R\cap A_q = \emptyset$ then $\phi(R)\cap\phi(A_q) = \phi(R)\cap
C_q^{\hat w} = \emptyset$, and $\phi(\Delta_q^{\hat w})$ is a pruning
disk which is disjoint from~$R$ but contains the highest points
$\{r_1,r_2,r_3,r_4\}$ of the orbits $P^w_q$.

Suppose, then, that $R\cap A_q\not=\emptyset$ and $R\cap
C_q\not=\emptyset$. It is required to prove that $R\ge P^w_q$. By
Lemma~\ref{lem:linking}, the linking number~$L_1$ is strictly smaller
than the linking numbers~$L_2$, $L_3$, and $L_4$. If~$L_1$ is not
divisible by~$n+k+3$, then an identical argument to that in the proof
of Theorem~\ref{thm:order} shows that~$(p_1\,;\,F)$ is unremovable
in~$(D^2,R)$, which establishes the result. On the other hand, if
$L_1$ is divisible by $n+k+3$ then $L_2$, $L_3$, and $L_4$ are
not. Either one of these linking numbers is distinct from the other
two, or all three are equal, and in either case the proof proceeds
identically to that of Theorem~\ref{thm:order}.

\end{proof}

\subsection{The invariant $r^w$}
The next result completes the proof of Theorem~\ref{thm:main} for
even decorations.
\begin{thm}
\label{thm:rw}
Let~$R$ be any horseshoe periodic orbit, and $r^w(R)\in(0,q_w]\cap\Q$
  be as given in Definition~\ref{defn:rw}. Then for $q\in\Q^w$
\begin{eqnarray*}
0<q<r^w(R)&\implies& R\not\ge P_q^w\qquad\text{and}\\
r^w(R)<q<q_w&\implies& R\ge P_q^w.
\end{eqnarray*}
In particular, $r^w$ is a braid type invariant.
\end{thm}
\begin{proof}  Let $q=m/n\in\Q^w$. Suppose first that
  $0<q<r^w(R)$: it is required to show that $R\not\ge P_q^w$. Without
  loss of generality, assume that $n/2$ is greater than the period
  of~$R$, so that Theorem~\ref{thm:forcing} applies (if not, replace
  $q$ with $q'\in(q,r^w(R))\cap\Q^w$ having sufficiently large
  denominator. Then $P_q^w\ge P_{q'}^w$ by Theorem~\ref{thm:order}, so
  showing that $R\not\ge P_{q'}^w$ shows that $R\not\ge P_q^w$).

Assume for a contradiction that $R\ge P_q^w$, so that $R\cap
A_q\not=\emptyset$ and $R\cap C_q\not=\emptyset$ by
Theorem~\ref{thm:forcing}. 

Let $x=b\cdot f\in R\cap C_q$. By Lemma~\ref{lem:gamma},
$f\succ\rib{c_q0w0}$ and $\sigma(b)\succ\rib{\hat{w}0c_q1}$. If
$\sigma(b)=\hat{w}0b'$ for some $b'\in\{0,1\}^\N$ then $b\cdot
f=b'0w\opti\cdot f$, where $f\succ\rib{c_q0w0}$ (so $q(f)\le q$) and
$b'\succ\rib{c_q1\hat{w}0}$ (so $q(b')\le q$). Thus there is a word
$010 w\opti10$ in $\overline{c_R}$ which forces $\lambda^w(R)\le q$.

On the other hand, if $\sigma(b)$ doesn't begin $\hat{w}0\ldots$ then
$\sigma(b)=\tilde{\hat{v}}b'$, where $\hat{v}$ is a non-empty even
initial subword of $\hat{w}0$. Then $b\cdot f=b'\breve{v}\opti\cdot f$
where $v$ is a non-empty even final subword of~$0w$, and
$f\succ\rib{c_q0w0}$ so that $q(f)\le q$. Hence $\mu^w(R)\le q$.

To summarize: the existence of a point of $R\cap C_q$ means that either
$\lambda^w(R)\le q$ or $\mu^w(R)\le q$.

Now let $x=b\cdot f\in R\cap A_q$ (to simplify the notation the same
symbols $x$, $b$, and~$f$ are used although~$x$ must be a different
point of~$R$). Thus $\phi(x)=f\cdot b\in\widehat{R}\cap
C_q^{\hat{w}}$, so that either $\lambda^{\hat{w}}(\widehat{R})\le q$
or $\mu^{\hat{w}}(\widehat{R})\le q$. Thus either $\lambda^w(R)\le q$
or $\nu^w(R)\le q$ by Remark~\ref{remark:rreverse}.

 So the fact that both $R\cap C_q$ and $R\cap A_q$ are non-empty means
 that either $\lambda^w(R)\le q$, or both $\mu^w(R)$ and $\nu^w(R)$
 are less than or equal to~$q$. This is precisely to say that
 $r^w(R)\le q$, which is the required contradiction.

For the converse, suppose that $r^w(R)<q<q_w$: it is required to show
that $R\ge P^w_q$. Once again, assume without loss of generality that
$n/2$ is greater than the period of~$R$ (if not, replace~$q$ with an
element of $(r^w(R),q)$ having sufficiently large denominator).

Since $q>r^w(R)=\min(\lambda^w(R),\max(\mu^w(R),\nu^w(R)))$, at least
one of the following two possibilities holds: that $q>\lambda^w(R)$,
or that $q$ is greater than both $\mu^w(R)$ and $\nu^w(R)$. It will be
shown that in either case~$R$ contains points of both~$A_q$ and~$C_q$,
so that $R\ge P^w_q$ by Theorem~\ref{thm:forcing} as required.

\begin{enumerate}[a)]
\item Suppose that $q>\lambda^w(R)$. Thus there is a point $x=b\cdot
  f$ of~$R$ with $q(f)<q$ and $b=\opti\hat{w}\opti b'$, where
  $q(b')<q$. In particular, $f\succ\rib{c_q0w0}$ and
  $\sigma(b)=\hat{w}\opti b'\succ\rib{\hat{w}0c_q1}$, so $x\in C_q$ by
  Lemma~\ref{lem:gamma}. Similarly,
  $q>\lambda^{\hat{w}}(\widehat{R})=\lambda^w(R)$, so there is a point
  of~$\widehat{R}$ in $C_q^{\hat{w}}$, and hence a point of~$R$ in~$A_q$.
\item Suppose that $q>\mu^w(R)$ and $q>\nu^w(R)$. The fact that
  $q>\mu^w(R)$ means that there is a point $x=b\cdot f$ of~$R$ with
  $q(f)<q$ and~$b=\opti\hat{\breve{v}}\ldots$, where~$v$ is a
  non-empty even final subword of $0w$. Thus $f\succ\rib{c_q0w0}$ and
  $\sigma(b)\succ\rib{\hat{w}0c_q1}$, so $x\in C_q$ by
  Lemma~\ref{lem:gamma}. Similarly,
  $q>\nu^w(R)=\mu^{\hat{w}}(\widehat{R})$ means that there is a point
  of~$\widehat{R}$ in $C_q^{\hat{w}}$, and hence a point of~$R$ in~$A_q$.
\end{enumerate}
\end{proof}
\subsection{The case of odd decoration}
\label{sec:odd}
The case where~$w$ is an odd decoration works similarly, the only
substantial changes being the exchange of the r\^oles of $C_q$
and~$D_q$ (and similarly of $A_q$ and~$B_q$), and the replacement of
the inclusions of Lemma~\ref{lem:contained} with ``approximate
inclusions''.

\begin{figure}[htbp]
\lab{1}{p_1}{r}
\lab{3}{\gamma}{r}
\lab{4}{S}{l}
\lab{5}{\delta}{l}
\lab{6}{p_3}{r}
\lab{7}{p_4}{r}
\lab{9}{p_2}{b}
\begin{center}
\pichere{0.4}{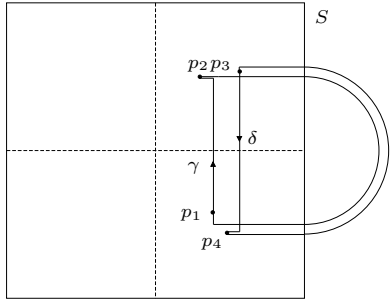}
\end{center}
\caption{The disks $C_q$ and $D_q$ when $w$ is odd}
\label{fig:alldiskodd}
\end{figure}

Define the points~$p_i$ and the arcs~$\gamma$, $\delta$, $\alpha$, and
$\beta$ exactly as in Section~\ref{sec:itarc}: observe that
now~$\Index(p_i,F^{n+k+3})=(-1)^{i+1}$ rather than $(-1)^i$ as
before. \mbox{Lemmas~\ref{lem:gamma}\,--\,\ref{lem:beta}} (describing the
disks defined by these arcs) and Lemma~\ref{lem:linking} and
Corollary~\ref{cor:linkcoincidence} (describing the linking numbers of
the orbits~$P_{i,q}$ about a periodic orbit~$R$) are
unchanged. However Lemma~\ref{lem:contained} is false, as can clearly
be seen from Figure~\ref{fig:alldiskodd}, which is a schematic
depiction of the points $p_i$ (with their correct relative horizontal
and vertical ordering), the arcs $\gamma$ and $\delta$, and the images
of these arcs under $F^{n+k+3}$. What is true, though, is that the
disk $C_q$ is {\em approximately} contained in $D_q$ in the sense of
Definition~\ref{defn:approx}.

\noindent\textbf{\em Lemma~\ref{lem:contained} (odd version).}\,\,---\,
\emph{For any odd decoration~$w$ and any~$q\in(0,q_w)$, $C_q\subseteq_q D_q$
and $A_q\subseteq_q B_q$.
}
\begin{proof}
Let~$q=m/n$. Let $x=\overline{c}$ be a point of a period~$N=|c|<n/2$
orbit~$R$ of~$F$. It suffices to show that if~$x\in\Int(C_q)$ then
$x\in\Int(D_q)$. The result will then follow since
$A_q=\phi(C^{\hat{w}}_q)$ and $B_q=\phi(D^{\hat{w}}_q)$. The proof
works in exactly the same way as the final part of the proof of
Lemma~\ref{lem:enlarge}. 

\end{proof}

Only a minor modification is needed to Lemma~\ref{lem:Cq}:

\noindent\textbf{\em Lemma~\ref{lem:Cq} (odd version).}\,\,---\,
\emph{Let $q,q'\in(0,q_w)\cap\Q$ with $q'<q$. Let~$R$ be the periodic orbit
of code~$c_{q'}0w0$. Then $|R\cap A_q|=|R\cap C_q|=0$, and $|R\cap
B_q|=|R\cap D_q|=1$. 
}

The proof works in exactly the same way, the only small difference
being the need to show that points of~$R$ other than the rightmost
point lie in {\em neither} $C_q$ nor $D_q$: this uses the fact that
if~$w$ is odd then $q(10w010\ldots)=q(10w110\ldots)$.

Theorem~\ref{thm:order} is true as stated for odd decorations, and the
proof is identical except that~$L_{3,q}(R)$ plays the r\^ole
of~$L_{1,q}(R)$, since the revised Lemma~\ref{lem:Cq} gives that
\[L_{1,q}(R)=L_{2,q}(R)=L_{4,q}(R)=L_{3,q}(R)-1.\]

The pruning disk which contains all of the points~$p_i$ is different
from that used in the case of even decoration.

\noindent\textbf{\em Lemma~\ref{lem:enlarge} (odd version).}\,\,---\,
\emph{Let $v_0$ and $v_1$ be the points $\li{}\,0c_q1w\opti\cdot
\rib{c_q0w10}$ if $w$ ends with the word $01^{2i}$ for some $i\ge 0$,
and the points $\li{}\,0c_q1w\opti\cdot
\rib{c_q0w110}$ otherwise. Then $v_0$ and $v_1$ are the vertices of a
pruning disk $\Delta_q$ which contains~$\{p_1,p_2,p_3,p_4\}$. Moreover
$\Delta_q \subseteq_q D_q$.
}

The only change to Theorem~\ref{thm:forcing} is the replacement of
$A_q$ and $C_q$ with $B_q$ and $D_q$. The proof works identically.

\noindent\textbf{\em Theorem~\ref{thm:forcing} (odd
  version).}\,\,---\, 
\emph{Let~$R$ be a period~$N$ horseshoe orbit, and suppose that
  $q=m/n\in\Q^w$ is such that $n>2N$. Then $R\ge P^w_q$ if and only if
  $R\cap B_q\not=\emptyset$ and $R\cap D_q\not=\emptyset$.  }

Finally, the statement of Theorem~\ref{thm:rw} is unchanged, and its
proof works in just the same way. Using~$B_q$ and $D_q$ in place
of~$A_q$ and~$C_q$ is exactly what is required to compensate for the
changes introduced because~$w$ is odd, and because $\lp{w}=1w$ and
$w^+=w1$, rather than $0w$ and $w0$ as in the even case.

\section{Examples and applications}
\label{sec:examples}
\subsection{Generalities}
The following straightforward lemma will be useful in this section.
\begin{lem}
\label{lem:gen}
Let~$w$ be a lone decoration. 
\begin{enumerate}[a)]
\item If $R$ is a horseshoe periodic orbit of decoration~$w$, then
$r^w(R)=q(R)$.
\item If~$R$ is any horseshoe periodic orbit and $r^w(R)\not=q_w$, then
$r^w(R)\ge q(R)$.
\item If~$R$ is any horseshoe periodic orbit, then
$r^w(R)=r^{\hat w}(\widehat{R})$.
\end{enumerate}
\end{lem}

\begin{proof}
a) is immediate from Theorem~\ref{thm:order}, Theorem~\ref{thm:rw},
and Lemma~\ref{lem:heightfacts}~c), while~c) is immediate from
Remark~\ref{remark:rreverse}.

For~b), suppose that $r^w(R)<q(R)$. Pick
$q_1,q_2\in(r^w(R),\min(q(R),q_w))\cap\Q^w$ with $q_1<q_2$. Since
$q_2<q(R)$, Lemma~\ref{lem:heightfacts}~b) gives $P^\ast_{q_2}\ge
R$. Since $r^w(R)<q_1$, Theorem~\ref{thm:rw} gives $R\ge
P^w_{q_1}$. Hence $P^\ast_{q_2}\ge P^w_{q_1}$, and
Lemma~\ref{lem:heightfacts}~c) gives $q_2=q(P^\ast_{q_2}) \le q(P^w_{q_1})
= q_1$, which is the required contradiction.
\end{proof}
\subsection{Star decorations}
\label{sec:star}
For each rational $m/n\in(0,1/2)$, consider the ``star'' decoration
$w_{m/n}$ which is defined by removing the initial symbols $10$ and
the final symbols $01$ from the word $c_{m/n}$ of
Definition~\ref{defn:cmn}. Thus, using the notation of
Lemma~\ref{lem:starlem},
\[w_{m/n}=0^{\ka_1-1}1^20^{\ka_2}1^20^{\ka_3}1^2\ldots1^20^{\ka_{m-1}}1^2
0^{\ka_m-1}. \]

These decorations were considered in~\cite{Stars}: their name is due
to the fact that the train tracks for periodic orbits with decoration
$w_{m/n}$ are all star-shaped, with $n$ branches: the pseudo-Anosov
has an $n$-pronged singularity corresponding to the vertex of the
star, whose prongs are rotated by $m/n$.

The scope $q_{w_{m/n}}$ of the decoration~$w_{m/n}$ is $m/n$. In
order to simplify the notation, the periodic orbits $P^{w_{m/n}}_q$ of
decoration~$w_{m/n}$ and height $q\in(0,m/n)$ will be denoted
$P^{m/n}_q$ and the invariants $r^{w_{m/n}}$ will be denoted $r^{m/n}$
throughout this subsection. 

It is shown in Lemma~17 and Corollary~18 of~\cite{Stars} that each
$w_{m/n}$ is a lone decoration, and that the periodic orbits
$P^{m/n}_q$ all have pseudo-Anosov braid type. Thus for each
fixed~$m/n$ the set
\[\cD^{m/n}=\{\bt(P_q^{m/n},F)\,:\,q\in(0,m/n)\}\]
is totally ordered by forcing, with $P_{q'}^{m/n}\ge P_q^{m/n}$ if and
only if $q'\le q$ (Theorem~\ref{thm:order}). The aim in this section
is to use Theorem~\ref{thm:rw} to determine the forcing between braid
types in distinct families $\cD^{m/n}$: that is, for each $m/n$,
$m'/n'$, $q$, and $q'$, to determine whether or not $P^{m'/n'}_{q'}\ge
P^{m/n}_q$. This completes the proof of Theorem~15~d) of~\cite{Stars}.

Recall that the {\em rotation number} $\rho(R)\in(0,1/2]$ of a horseshoe
periodic orbit~$R$ (other than a fixed point) is its rotation number
in the annulus obtained by puncturing the disk at the fixed point of
code $1$. The {\em rotation interval} $\rotint(R)\subseteq(0,1/2]$ is
  the set of rotation numbers of periodic orbits forced
  by~$R$. Clearly $\rotint(R)$ is a braid type invariant, and if $R\ge
  S$ then $\rotint(R)\supseteq\rotint(S)$.

Lemma~17 of~\cite{Stars} states that $\rotint(P_q^{m/n})=[q,m/n]$.
Hence if $P_{q'}^{m'/n'}\ge P_q^{m/n}$ then $[q,m/n]\subseteq
[q',m'/n']$. The main result of this section is a near-converse to
this statement.

\begin{thm}
\label{thm:starforce}
Let $m/n, m'/n', q' \in (0,1/2)\cap\Q$ with $q'<m'/n'$. Then
\[r^{{m/n}}(P_{q'}^{m'/n'})=\left\{
\begin{array}{ll}
q' & \text{ if }q'<\frac{m}{n}\le\frac{m'}{n'} \\
\frac{m}{n} & \text{ otherwise.}
\end{array}
\right.\]
\end{thm}

That is, if $q'\not=q$ then $P^{m'/n'}_{q'}\ge P^{m/n}_q$ if and only
if $[q,m/n]\subseteq [q',m'/n']$. In the language of
Conjecture~\ref{conj:dec}, this means that $w_{m/n}\preceq w_{m'/n'}$
if and only if $m/n\le m'/n'$.

\begin{proof} 
If $m/n \le q'$ or $m/n>m'/n'$, then
$\rotint(P_{q'}^{m'/n'})=[q',m'/n']$ does not contain the rotation
interval $[q,m/n]$ of $P_q^{m/n}$ for any~$q$: that is,
$P_{q'}^{m'/n'}$ does not force any $P_q^{m/n}$, and hence
$r^{{m/n}}(P_{q'}^{m'/n'})=m/n$ as required.

If $m'/n'=m/n$, then $r^{{m/n}}(P_{q'}^{m'/n'})=q'$ as required by
Lemma~\ref{lem:gen}~a).

Suppose, then, that $q'<m/n<m'/n'$. Lemma~\ref{lem:gen}~b) gives that
$r^{{m/n}}(P_{q'}^{m'/n'})\ge q'$, so it suffices to show that
$r^{{m/n}}(P_{q'}^{m'/n'})\le q'$. To do this, it is enough to show
that $\mu^{w_{m/n}}(P_{q'}^{m'/n'})\le q'$ and
$\nu^{w_{m/n}}(P_{q'}^{m'/n'})\le q'$.

Let $R=P_{q'}^{m'/n'}$ with $c_R=c_{q'}0w_{m'/n'}0$, and write
\[w_{m'/n'}=0^{\ell_1-1}1^20^{\ell_2}1^20^{\ell_3}1^2\ldots 
1^20^{\ell_{m'-1}}1^20^{\ell_{m'}-1}.\] Similarly, write
\[w_{m/n}=0^{\kappa_1-1}1^20^{\kappa_2}1^20^{\kappa_3}1^2\ldots 
1^20^{\kappa_{m-1}}1^20^{\kappa_{m}-1}.\]

Consider the point of~$R$ with itinerary
$\sigma^{n'+1}(\overline{c_R})=\overline{0w_{m'/n'}0c_{q'}}=b\cdot f$,
so that \mbox{$f=\rib{0w_{m'/n'}0c_{q'}}$} and
$b=\rib{c_{q'}0w_{m'/n'}0}$ (the latter using the fact that
$w_{m'/n'}$ and $c_{q'}$ are palindromic). It will be shown that
$f=0\tilde v f'$ for some $f'\in\{0,1\}^\N$, where $v$ is a non-empty
even initial subword of $w_{m/n}^+$. Thus
\[
\sigma^{n'+|v|}(\overline{c_R})=b0\tilde v\cdot f',
\]
for $v$ a non-empty even initial subword of $w_{m/n}^+$, so
\mbox{$\nu^{w_{m/n}}(P_{q'}^{m'/n'})\le q(b)=q'$}. Then
\[\mu^{w_{m/n}}(P_{q'}^{m'/n'})=\nu^{\hat{w}_{m/n}}(\hat{P}_{q'}^{m'/n'})=\nu^{w_{m/n}}(P_{q'}^{m'/n'})
\le q',\]
so that $r^{{m/n}}(P_{q'}^{m'/n'})\le q'$ as required.

Thus it only remains to establish the claim, that
$f=\rib{0w_{m'/n'}0c_{q'}}$ is of the form $0\tilde v f'$, where
$v$ is a non-empty even initial subword of $w_{m/n}^+=w_{m/n}0$.

Observe first that the words
\begin{eqnarray*}
w_{m'/n'}0c_{q'}&=&0^{\ell_1-1}1^20^{\ell_2}1^2\ldots
1^20^{\ell_{m'-1}}1^2 0^{\ell_{m'}}10\ldots \quad\text{ and}\\
w_{m/n}0&=&0^{\ka_1-1}1^20^{\ka_2}1^2\ldots 1^20^{\ka_{m-1}}1^20^{\ka_m}
\end{eqnarray*}
must disagree before the shorter of their lengths. For
\begin{enumerate}[a)]
\item There is a subword of the form $010$ in $w_{m'/n'}0c_{q'}$ but no
such subword in $w_{m/n}0$. Hence if $|w_{m/n}0|\ge|w_{m'/n'}0c_{q'}|$
then the two words must disagree.
\item On the other hand, if $|w_{m/n}0|<|w_{m'/n'}0c_{q'}|$ and the
  words do not disagree, then $w_{m'/n'}0c_{q'}=w_{m/n}0\ldots$, and
  so $10w_{m'/n'}0c_{q'}=10w_{m/n}0\ldots$ and (since $10w_{m/n}0$ is
  an odd word)
\[10w_{m'/n'}0c_{q'} \succ 10w_{m/n}010^\infty.\]
Taking the height of both sides gives $m'/n' \le m/n$, a contradiction.
\end{enumerate}

Since $m'/n'=q(10w_{m'/n'}0c_{q'}\ldots)>q(10w_{m/n}010\ldots)=m/n$,
the word $w_{m'/n'}0c_{q'}$ is greater than $w_{m/n}0$ in the unimodal
order. Let~$u$ be the longest initial word on which they agree,
and~$v$ be the length~$|u|+1$ initial subword of $w_{m/n}0$. Then
either $u$ is even and $v=u0$, or $u$ is odd and $v=u1$. In either
case, $v$ is a non-trivial even initial subword of
$w_{m/n}^+=w_{m/n}0$ and $0\tilde v$ is an initial subword of
$0w_{m'/n'}0c_{q'}$ as required.
\end{proof}

\subsection{Decorations of the form $1^{2i+1}$}
\label{sec:111}
For each~$i\ge 0$, consider the decoration $w_i=1^{2i+1}$, with scope
$q_{w_i}=q\left(\rib{101^{2i+1}0}\right)=1/2$. It will be shown
that~$w_i$ is lone for all~$i$, yielding corresponding braid type
invariants $r^{w_i}$ by Theorem~\ref{thm:main}. These invariants will
then be used to show that $P^{w_i}_q$ has pseudo-Anosov braid type for
all~$i$ and all~$q\in(0,1/2)$.

In order to simplify the notation, the periodic orbits $P_q^{w_i}$
will be denoted $P_q^i$ and the invariants $r^{w_i}$ will be denoted
$r^i$ throughout this subsection. Similarly, $\lambda^{w_i}$, $\mu^{w_i}$
and $\nu^{w_i}$ will be abbreviated to $\lambda^i$, $\mu^i$, and
$\nu^i$. 

\begin{remark}
  Decorations of the form $1^{2i}$ ($i\ge 0$) are also lone: $1^{2i}$
  is the decoration $w_{(i+1)/(2i+3)}$ of Section~\ref{sec:star}.
\end{remark}

The proof that the decorations~$w_i$ are lone is by induction on~$i$:
the fact that $w_i$ is lone will be established by using the fact that
$w_{i-1}$ is lone, and hence that $r^{i-1}$ is a braid type
invariant. The proof will also make use of the following theorem and
lemma, which can be found in~\cite{Ha}, Theorem~\ref{thm:ri} appearing
as Theorems~3.11 and~3.15, and Lemma~\ref{lem:ht12} appearing as
Lemma~3.3.

\begin{thm}
\label{thm:ri}
The rotation interval of any horseshoe periodic orbit~$R$ is of the
form
\[\rotint(R)=[q(R),\operatorname{rhe}(R)]\cap\Q,\]
where $\operatorname{rhe}(R)\in[q(R),1/2]\cap\Q$ is equal to~$1/2$ if
and only if $\overline{c_R}$ contains either the word~$01010$ or the
word~$01^{2m+1}0$ for some $m\ge 1$.
\end{thm}

\begin{lem}
\label{lem:ht12}
Suppose that $c\in\{0,1\}^\N$ has $q(c)=1/2$. Then $c=0\ldots$, or
\mbox{$c=11\ldots$}, or $c=101^{2k-1}0\ldots$ for some $k\ge1$.
\end{lem}

The following lemma will also be used in the proof.
\begin{lem}
\label{lem:ri12}
Let $i\ge 0$, and let~$R$ be a horseshoe periodic orbit with
$r^i(R)<1/2$. Then $\overline{c_R}$ contains a word of the form
$01^{2k+1}0$ for some $k\le i+2$.
\end{lem}
\begin{proof}
(Note that it is not assumed in this proof that $r^i$ is a braid type
  invariant.) 

Since $r^i(R)=\min(\lambda^i(R),\max(\mu^i(R),\nu^i(R)))<1/2$, either
$\lambda^i(R)<1/2$ or \mbox{$\nu^i(R)<1/2$}.

If $\lambda^i(R)<1/2$, then there is a word $v=\opti\, w^i\,\opti$ such
that some shift of $\overline{c_R}$ is of the form $bv\cdot f$, where
$q(b)<1/2$ and $q(f)<1/2$ (i.e. $b=10\ldots$ and $f=10\ldots$). Thus
$bv\cdot f=\ldots01\opti 1^{2i+1}\opti 10\ldots$, which contains a
block of $1$s of odd length either $2i+1$, $2i+3$, or $2i+5$ as
required.

Similarly, if $\nu^i(R)<1/2$ then there is some non-empty even final
subword $v=1^{2j}$ \mbox{($1\le j\le i+1$)} of $w^{i+}=1^{2i+2}$ with
the property that some shift of $\overline{c_R}$ is of the form
\mbox{$b\opti\tilde{v}\cdot f = b\opti 1^{2j-1}0\cdot f$}, where
$q(b)<1/2$ and so $b=10\ldots$. Hence $\overline{c_R}$ contains a
block of $1$s of odd length either $2j-1$ or $2j+1\le 2i+3$ as
required.
\end{proof}

\begin{thm}
\label{thm:lone}
For each~$i\ge 0$, the decoration~$w_i=1^{2i+1}$ is lone.
\end{thm}
\begin{proof}
The proof is by induction on~$i$. For~$i=0$, it is required to show
that for each $q=m/n\in(0,1/2)\cap\Q$, the four horseshoe periodic
orbits~$P^{w_0}_q$ of codes $c_q\opti 1\opti$ are the only horseshoe
periodic orbits of their braid types. Now the only other horseshoe
periodic orbits of height~$q$ and period~$n+4$ are those of
decoration~$0$, i.e. those with codes $c_q\opti 0\opti$ (since
$q_0=1/3$, there is in fact nothing to prove if $q>1/3$). However
$\rotint(P^{w_0}_q)=[q,1/2]$ by Theorem~\ref{thm:ri} (since
$\overline{c_q\opti1\opti}$ contains the word $01\opti1\opti10$),
while the periodic orbits of codes $c_q\opti 0\opti$ have rotation
intervals with right hand endpoints less than $1/2$, since
$\overline{c_q\opti0\opti}$ contains neither of the words of
Theorem~\ref{thm:ri}. (In fact these are periodic orbits with star
decoration $w_{1/3}$, and hence their rotation intervals have right
hand endpoint $1/3$.)

Now let~$i>0$. By the inductive hypothesis, $r^{i-1}$ is a braid type
invariant. 

Lemma~\ref{lem:ri12} gives $r^{i-1}(P^i_q)=1/2$. For writing
$R=P^i_q$, with
\[c_R=c_q1w_i1=10^{\kappa_1}1^2\ldots 1^20^{\kappa_m}1^{2i+4},\]
$\overline{c_R}$ contains only one block of $1$s of odd length, and
that block has length $2i+5$.

Thus $P^i_q$ is a period~$n+2i+4$ orbit with rotation
interval~$[q,1/2]$ and $r^{i-1}(P^i_q)=1/2$. It will be shown that any
horseshoe orbit with these properties must be $P^i_q$, which will
complete the proof.

Suppose then that~$R$ is a period~$n+2i+4$ periodic orbit with height
$q(R)=q$, with $r^{i-1}(R)=1/2$, and for which (using
Theorem~\ref{thm:ri}) $\overline{c_R}$ contains either the word
$01010$ or the word $01^{2m+1}0$ for some $m\ge 1$. Take 
\[c_R=c_q1w'1 = 10^{\kappa_1}1^2\ldots1^20^{\kappa_m}1w'1\]
for some decoration $w'$ of length $2i+1$ distinct from $1^{2i+1}$. A
contradiction will be derived.

Consider the first occurrence of one of the words $01010$ or
$01^{2m+1}0$ ($m\ge 1$) in $c_R^\infty$. Since $w'\not=1^{2i+1}$, this
word is either $01010$ or $01^{2m+1}0$ for $1\le m\le i$. Thus there
is a corresponding shift of $\overline{c_R}$ of the form $\ldots
01\opti 1^{2m-1}0\cdot f$ where $1\le m\le i$; or, in other words, of
the form $b\opti\tilde{v}\cdot f$, where $v=1^{2m}$ is a non-trivial
even initial subword of $w^{(i-1)+}=1^{2i}$.

Because we chose the {\em first} occurrence of one of the words $01010$
or $01^{2m+1}0$ in $c_R^\infty$, $b=10\ldots$ cannot be of the form
$101^{2k-1}0\ldots$ for any $k\ge 1$. Thus $q(b)<1/2$ by
Lemma~\ref{lem:ht12}, and hence $\nu^{i-1}(R)<1/2$.

Now $\mu^{i-1}(R)=\nu^{\hat{w}_{i-1}}(\hat{R}) =
\nu^{i-1}(\hat{R})$. Since $q(\hat{R})=q(R)$ (Lemma~3.8 of~\cite{Ha})
and $\overline{c_{\hat{R}}}$ contains one of the words $01010$ or
$01^{2m+1}0$, this gives $\mu^{i-1}(R)<1/2$. Hence $r^{i-1}(R)<1/2$,
which is the required contradiction.

\end{proof}

Thus $r^i$ is a braid type invariant for all~$i\ge 0$. These
invariants will now be used to prove that $P^i_q$ has pseudo-Anosov
braid type for all $i\ge 0$ and $q\in(0,1/2)$. The proof will make use
of a more general result (Theorem~\ref{thm:pa} below) for showing that
horseshoe periodic orbits have pseudo-Anosov braid type.

The following result (Lemma~3.4 of~\cite{Ha}) will be required:
\begin{lem}
\label{lem:heightchar}
Let $q\in(0,1/2)\cap\Q$. Let $w_q$ be as in
Section~\ref{sec:star}. Then, for any $c\in\{0,1\}^\N$,
\[q(c)=q \iff (10w_q1)^\infty \preceq c \preceq
10w_q0(11w_q0)^\infty.\] 
In particular, if $q(c)=q$ then $c=10w_q110\ldots$ or $c=10w_q01\ldots$.
\end{lem}

The first step is to give a lower bound on the period of horseshoe
orbits~$R$ for which $r^i(R)<r^{i-1}(R)$.

\begin{lem}
\label{lem:periodbound}
Let $q=m/n\in(0,1/2)$, and let~$R$ be a horseshoe periodic orbit such
that $r^i(R)=q<r^{i-1}(R)$. Then~$R$ has period at least $n+2i+4$.
\end{lem}
\begin{proof} 
\begin{enumerate}[a)]
\item Suppose first that $\lambda^i(R)=q$. It follows from
Definitions~\ref{defn:forwardbackward} that for some~$s$
\[\sigma^s(\overline{c_R})=b\opti 1^{2i+1}\opti\cdot f,\]
where $q(b)\le q$ and $q(f)\le q$, with equality in one of the two
cases. Suppose that $q(f)=q$: the case where $q(b)=q$ works
identically. Hence (Lemma~\ref{lem:heightchar}) either
$f=10w_q110\ldots$ or $f=10w_q01\ldots$. In the former case,
$\overline{c_R}$ contains both the word $01\opti 1^{2i+1}\opti 10$
(which has length $2i+7$ and contains only isolated $0$s and blocks of
$1$s of odd length) and the word $0w_q110$ (which has length $n+1$ and
contains only blocks of $1$s of even length). Thus~$R$ has period at
least $(2i+7)+(n+1)-2 > n+2i+4$ as required. In the latter case,
$\overline{c_R}$ again contains the word $01\opti 1^{2i+1}\opti 10$,
and also contains the word $0w_q0$ (which has length $n-1$ and
contains only blocks of $1$s of even length). Thus~$R$ has period at
least $(2i+7)+(n-1)-2 = n+2i+4$ as required. (Note that if~$R$ has
period $n+2i+4$, then $\sigma^s(\overline{c_R})=10w_q01\opti 1^{2i+1}
\opti=c_q\opti 1^{2i+1}\opti$, i.e. $R=P^i_q$.)

\item
Suppose, then, that $\lambda^i(R)\not=q$. Since $r^i(R)=q$, it follows
from Definitions~\ref{defn:forwardbackward} that $\lambda^i(R)>q$, and
that $\mu^i(R)\le q$ and $\nu^i(R)\le q$, with equality in one of the two
cases. Suppose that $\mu_i(R)=q$: the case where $\nu_i(R)=q$ works
identically. Hence, by Definitions~\ref{defn:forwardbackward}, there
is some~$s$ such that
\[\sigma^s(\overline{c_R}) = b01^{2k-1}\opti\cdot f,\]
where $1\le k\le i+1$ and $q(f)=q$ (here $01^{2k-1}=\breve{v}$, where
$v=1^{2k}$ is a non-empty even final subword of $\lp w^i = 1^{2i+2}$).
\begin{enumerate}[i)]
\item If $k=i+1$, then $\overline{c_R}$ contains the word
  $01^{2i+1}\opti 10$ (of length $2i+5$, with only isolated $0$s and
  blocks of $1$s of odd length). It also contains either the word
  $0w_q110$ (length $n+1$) or $0w_q0$ (length $n-1$), which contain
  only blocks of $1$s of even length. So if~$R$ has period less than
  $n+2i+4$, $\overline{c_R}$ must contain the words $01^{2i+1}\opti
  10$ and $0w_q0$, and these words must overlap at either one or both
  of their endpoints: that is, $\sigma^s(\overline{c_R})$ is either
  $\overline{10w_q01^{2i+1}\opti}$, or
  $\overline{10w_q001^{2i+1}\opti}$, or
  $\overline{100w_q01^{2i+1}\opti}$. In the first case, $R=P_q^{i-1}$,
  and hence $r^{i-1}(R)=q$, contradicting the hypothesis that
  $q<r^{i-1}(R)$. In the second case, $f=10w_q00\ldots$ and in the
  third case, $f=100w_q0\ldots$, each contradicting $q(f)=q$ (since by
  Lemma~\ref{lem:heightchar}, if $q(f)=q$ then the number of $1$s in
  the first $n+1$ symbols of~$f$ is either $2m$ or $2m+1$).
\item If $k\le i$, then $\mu^{i-1}(R)\le q$ (since $1^{2k}$ is a
  non-empty even final subword of $\lp w^{i-1}$). Since
  $r^{i-1}(R)>q$, it follows that $\nu^{i-1}(R)>q$.

Let $\nu^i(R)=r\le q$. Then $\overline{c_R}$ contains the word
$w_r01\opti 1^{2i+1}0$ (the block of $1$s here can't be shorter, since
$\nu^{i-1}(R)>r$). In particular, it contains the word $01\opti
1^{2i+1}0$ (which has length $2i+5$ and contains only isolated $0$s
and blocks of $1$s of odd length). Since $\mu^i(R)=q$,
$\overline{c_R}$ also contains either the word $0w_q110$ (length
$n+1$) or $0w_q0$ (length $n-1$). So if~$R$ has period less than
$n+2i+4$, $\overline{c_R}$ must contain the words $01\opti 1^{2i+1}0$
and $0w_q0$, and these words must overlap at either one or both of
their endpoints: that is, $\sigma^s(\overline{c_R})$ is either
$\overline{10w_q01\opti 1^{2i}}$, or $\overline{10w_q001\opti
  1^{2i}}$, or $\overline{100w_q01\opti 1^{2i}}$. As before, in the
first case $R=P_q^{i-1}$, contradicting $r^{i-1}(R)<q$, while the
other two cases contradict $q(f)=q$. 
\end{enumerate}

\end{enumerate}

\end{proof}

Note that there is no restriction on the decoration of the horseshoe
periodic orbit~$R$ in the following result, which thus provides a
general test for pseudo-Anosov braid type of horseshoe periodic orbits.

\begin{thm}
\label{thm:pa}
Let~$i\ge 1$ and let~$R$ be a period~$N$ horseshoe orbit with
\mbox{$r^i(R)=m/n<r^{i-1}(R)$}. Let~$d$ be the largest divisor of~$N$
other than~$N$ itself. If $d<n+2i+4$, then~$R$ has pseudo-Anosov braid
type.
\end{thm}
\begin{proof}
The fact that $r^i(R)<r^{i-1}(R)\le 1/2$ means that $R$ forces
infinitely many periodic orbits of decoration $w^i$, so cannot be of
finite order braid type.

If~$R$ had reducible braid type, then it would force the braid type of
the outermost component in its Nielsen-Thurston canonical
representative~$g$: in particular, this is the braid type of some
horseshoe periodic orbit~$S$. Since $R\ge S$, it follows that
$r^{i-1}(S)\ge r^{i-1}(R)$. Now~$R\ge P^i_q$ for all $q>r^i(R)$ with
$q\in\Q^{w_i}$, so $g$ has periodic orbits of each of these braid
types in its outermost component, and hence $S\ge P^i_q$ for all
$q>r^i(R)$ with $q\in \Q^{w_i}$. So $r^i(S)=r^i(R)<r^{i-1}(R)\le
r^{i-1}(S)$, and hence~$S$ has period at least~$n+2i+4$ by
Lemma~\ref{lem:periodbound}. This contradicts the fact that the period
of~$S$ is at most~$d$, which is less than $n+2i+4$.
\end{proof}

\begin{cor}
\label{cor:pa}
$P^i_q$ has pseudo-Anosov braid type for all $i\ge 0$ and all $q\in(0,1/2)$.
\end{cor}
\begin{proof}
Let $q=m/n$. Suppose first that $i>0$. Then $r^{i-1}(P^i_q)=1/2$,
as established in the proof of Theorem~\ref{thm:lone}; and
$r^i(P^i_q)=q$ by Lemma~\ref{lem:gen}~a). Since $P^i_q$ has period
$n+2i+4$, the result follows from Theorem~\ref{thm:pa}.

For the case~$i=0$, suppose for a contradiction that~$R$ has reducible
braid type. As in the proof of Theorem~\ref{thm:pa}, let~$S$ be a
horseshoe periodic orbit whose braid type is that of the outermost
component in the Nielsen-Thurston canonical representative of the
braid type of~$R$. Then $\rotint(S)=\rotint(R)$, so in particular
$q(S)=q(R)=m/n$, and hence~$S$ has period at least~$n$. Since~$R$ has
period~$n+4$, the period of~$S$ is at most $n/2 + 2$: thus $n\le 4$. A
direct check verifies that the orbits $P^0_{1/3}$ and $P^0_{1/4}$ have
pseudo-Anosov braid type.
\end{proof}

It follows from Theorem~\ref{thm:main} that $P^i_q\ge P^i_{q'}$ for
all~$i$ and all $q,q'\in(0,1/2)$ with $q\le q'$. The forcing between
families with different~$i$ can also be determined easily using the
invariants~$r^i$. The next result says that $P^i_q$ forces none of the
$P^{j}_{q'}$ with $j<i$, while if $j\ge i$ it forces all those with
$q'>q$: in the language of Conjecture~\ref{conj:dec}, this means that
$w^{j}\preceq w^i$ if and only if $j\ge i$.

\begin{thm}
\label{thm:interwiforcing}
Let $i$ and $j$ be non-negative integers, and
$q\in(0,1/2)\cap\Q$. Then 
\[r^j(P^i_q) = \left\{
\begin{array}{ll}
q & \text{ if }j\ge i\\
\frac12 & \text{ if }j< i.
\end{array}
\right.
\]
\end{thm}
\begin{proof}
Let $R=P^i_q$ with $c_R=c_q1^{2i+3}$. Then the only words of the form
$01^{2k+1}0$ in $\overline{c_R}$ have $k=i+2$. It is therefore
immediate from Lemma~\ref{lem:ri12} that $r^j(R)=1/2$ for $j<i$.

$r^i(R)=q$ by Lemma~\ref{lem:gen}~a), so suppose that
$j>i$. $r^j(R)\ge q$ by Lemma~\ref{lem:gen}~b), so it is only
necessary to show that $r^j(R)\le q$. Now $\overline{c_R}=\ldots
01^{2i+3}1\cdot\rib{c_q1^{2i+3}}$, and $01^{2i+3}1$ is a word of the
form $\breve{v}\opti$, where $v=1^{2i+4}$ is a non-empty even final
subword of $\lp{w}^j=1^{2j+2}$. Hence $\mu^j(R)\le
q(\rib{c_q1^{2i+3}})=q$. Similarly $\nu^j(R)\le q$, giving $r^j(R)\le
q$ as required.
\end{proof}

\begin{cor}
\label{cor:ridecreases}
Let $R$ be a period~$N$ horseshoe orbit. Then $(r^i(R))$ is a
decreasing sequence, with $r^i(R)=r^{i'}(R)$ if $i,i'\ge
\lfloor(N-7)/2\rfloor$.
\end{cor}
\begin{proof}
Let $q>r^i(R)$ and pick $q'\in(r^i(R),q)$. Then $R\ge P^i_{q'}$
by Theorem~\ref{thm:main}, and $P^i_{q'}\ge P^{i+1}_q$ by
Theorem~\ref{thm:interwiforcing}. That is, $R\ge P^{i+1}_q$ for all
$q>r^i(R)$, and so $r^{i+1}(R)\le r^i(R)$ as required.

That the sequence $(r^i(R))$ stabilises after
$i=\lfloor(N-7)/2\rfloor$ is immediate from
Lemma~\ref{lem:periodbound}.
\end{proof}

\subsection{Topological entropy bounds}
\label{sec:entropy}

Recall that the {\em topological entropy} $h(\beta)$ of a braid
type~$\beta$ is the minimum topological entropy of
orientation-preserving homeomorphisms of the disk having a periodic
orbit of braid type~$\beta$: it is realised by the Nielsen-Thurston
canonical representative of the braid type.

Let~$w$ be a lone decoration, and let $h^w(q)=h(P^w_q)$ denote the
topological entropy of the braid type of the periodic orbits with
height~$q$ and decoration~$w$ ($0<q<q_w$). It is clear that, for any
horseshoe periodic orbit~$R$, $h(R)\ge h^w(q)$ for all $q>r^w(R)$ and,
in particular, that $h(R)\ge \bh^w(r^w(R))$, where
\[\bh^w(q)=\lim_{q'\searrow q}h^w(q').\]

It is often possible to calculate $\bh^w(q)$ explicitly using train
track techniques, providing a convenient means to compute topological
entropy bounds. The approach for the decorations $w_i$ of
Section~\ref{sec:111} will be outlined in this section. A similar
calculation could in principle be carried out for the star decorations
of Section~\ref{sec:star}, using the explicit train track maps
described in~\cite{Stars}.

\begin{figure}[htbp]
\begin{center}
\pichere{0.9}{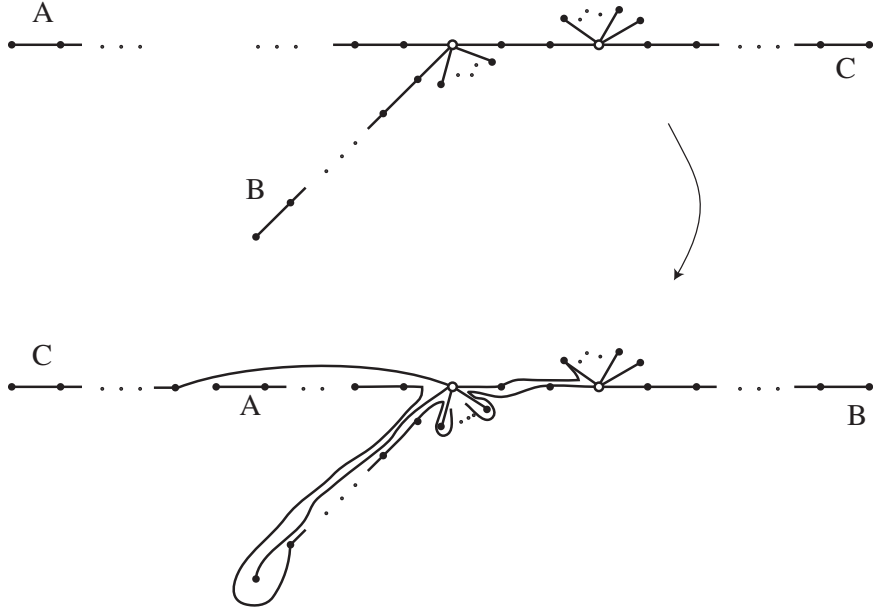}
\end{center}
\caption{The train track map for $P_q^{w_i}$}
\label{fig:tt}
\end{figure}

An explicit train track and train track map for the periodic
orbits $P^{w_i}_q$ is depicted in Figure~\ref{fig:tt}. Writing $q=m/n$,
the $n+2i+4$ points of the orbit are depicted with solid
circles. There are two valence $i+3$ vertices, depicted with unfilled
circles. The strings of edges denoted A, B, and C contain respectively
$n-2m+1$, $m$, and $m$ points of the orbit (the remaining $2i+3$
points comprising $i$ at valence~1 vertices around the left hand
valence $i+3$ vertex, $i+1$ at valence~1 vertices around the right
hand valence $i+3$ vertex, and $2$ between these two vertices).

A routine but long calculation using this train
track map shows that $h^{w_i}(m/n)$ is the logarithm of the largest
real root of the polynomial
\[H^i_{m/n}(x)=x^{n+1}g_i(x)+2x(x^2-1)(x^{2i+4}+1)f_{m/n}(x)-x^{2i+6}g_i(1/x),\]
where
\[g_i(x)=x^{2i+3}(x^3-x^2-x-1)-2 \quad\text{ and
}f_{m/n}(x)=\sum_{j=1}^{m-1}x^{\lfloor jn/m\rfloor}.\]
(In fact $H^i_{m/n}(x)$ is $(x^2-1)$ times the characteristic
polyomial of the transition matrix of the train track map.)

Using this result, the following theorem can be proved:
\begin{thm}
\label{thm:entropy}
Let $0<m/n<1/2$. The polynomial
\[\overline{H}^i_{m/n}(x)=(x^n-1)g_i(x)+2(x^2-1)(x^{2i+4}+1)(1+f_{m/n}(x))\]
has a single real root $\lambda^i_{m/n}$ in $x>1$, and
$\bh^{w_i}(m/n)=\log \lambda^i_{m/n}$.
\end{thm}
\begin{proof}
Let
\[K(x)=\frac{\overline{H}^i_{m/n}(x)}{x-1} = g_i(x)\sum_{j=0}^{n-1}x^j
+ 2(x+1)(x^{2i+4}+1)(1+f_{m/n}(x)).\]
Observe that $K(1)=4(2m-n)<0$ and that $K(2)>0$ since $g_i(2)\ge0$, so
that~$K$ has at least one root in~$(1,2)$. Now
\begin{eqnarray*}
K(x) &=& x^{2i+6}\sum_{j=0}^{n-1}x^j + 2(x+1)(x^{2i+4}+1)
\sum_{j=0}^{n-1}x^{\lfloor jn/m\rfloor} -
(x^{2i+5}+x^{2i+4}+x^{2i+3}+2)\sum_{j=0}^{n-1}x^j \\
&=& x^{n+2i+5}+2(x+1)(x^{2i+4}+1) \sum_{j=0}^{m-1}x^{\lfloor
  jn/m\rfloor} - x^{2i+5}-(x^{2i+4}+x^{2i+3}+2)\sum_{j=0}^{n-1}x^j.
\end{eqnarray*}
Since $n>2m$, every term of $2(x+1)(x^{2i+4}+1) \sum_{j=0}^{m-1}x^{\lfloor
  jn/m\rfloor}$ is cancelled by terms of
$(x^{2i+4}+x^{2i+3}+2)\sum_{j=0}^{n-1}x^j$, so that~$K(x)$ is of the
form
\[K(x)=x^{n+2i+5}-\sum_{j=0}^{n+2i+3}a_jx^j,\]
where all of the coefficients $a_j$ are non-negative. It follows that
all of its positive roots occur with positive derivative, so
that~$K(x)$ has a unique positive root, and hence
$\overline{H}^i_{m/n}$ has a unique root in $(1,\infty)$ as required.

Now $\bh^{w_i}(m/n)$ is the limit as $k\to\infty$ of the increasing
sequence $\mu_k=h^{w_i}\left(\frac{m^2k}{mnk-1}\right)$. Pick a
rational $m'/n'\in(0,1/2)$ which is greater than~$m/n$, and let
$\mu=h^{w_i}(m'/n')$, so that $\mu<\mu_k$ for all sufficiently
large~$k$. Restrict to such large~$k$, and work with values of~$x$ in
the interval~$[e^\mu,2]$. Then
\begin{eqnarray*}
H^i_{m^2k/(mnk-1)}(x) &=& x^{mnk}g_i(x)+2x(x^2-1)(x^{2i+4}+1)
\sum_{j=1}^{m^2k-1} x^{\left\lfloor
\frac{(m^2k-j)(mnk-1)}{m^2k}
\right\rfloor}
-x^{2i+6}g_i(1/x)\\
&=&
x^{mnk}\left(
g_i(x)+2(x^2-1)(x^{2i+4}+1)
\sum_{j=1}^{m^2k-1} x^{\left\lfloor
\frac{-j(mnk-1)}{m^2k}
\right\rfloor}
-x^{2i+6-mnk}g_i(1/x)
\right).
\end{eqnarray*}

Now
\[\sum_{j=1}^{m^2k-1} x^{\left\lfloor
\frac{-j(mnk-1)}{m^2k}
\right\rfloor} = 
\sum_{r=0}^{mk-1}x^{-rn}
\sum_{j=1}^{m-1}x^{\left\lfloor
\frac{-jn}{m}+\frac{rm+j}{m^2k}
\right\rfloor}
+\sum_{r=1}^{mk-1}x^{-rn},
\] 
and since $(rm+j)/m^2k<1/m$ for $r\le k-1$ and $j\le m-1$; and $(rm+j)/m^2k<1$ for
$r\le mk-1$ and $j\le m-1$, this gives
\begin{eqnarray*}
\sum_{j=1}^{m^2k-1} x^{\left\lfloor
\frac{-j(mnk-1)}{m^2k}
\right\rfloor} &=&
\sum_{r=0}^{k-1}x^{-rn}\sum_{j=1}^{m-1}x^{\lfloor-jn/m\rfloor} +
\sum_{r=1}^{mk-1}x^{-rn} + R_{k,m/n}(x)\\
&=&f_{m/n}(x)\sum_{r=1}^kx^{-rn}+\sum_{r=1}^{mk-1}x^{-rn} + R_{k,m/n}(x),
\end{eqnarray*}
where the remainder term $R_{k,m/n}(x)$ satisfies
\[0\le R_{k,m/n}(x)\le xf_{m/n}(x)\sum_{r=k+1}^{mk}x^{-rn} \quad\text{
  for all }x\in[e^\mu,2].\]

Thus
\[
\frac
{(x^n-1)H^i_{m^2k/(mnk-1)}(x)}
{x^{mnk}}
=\overline{H}^i_{m/n}(x)+S_{k,i,m/n}(x),
\]
where $S_{k,i,m/n}(x)\to 0$ as $k\to\infty$ uniformly for
$x\in[e^\mu,2]$. For each~$k$ this function has a zero at $e^{\mu_k}$,
and hence the unique zero in $(1,\infty)$ of $\overline{H}^i_{m/n}(x)$
is at $\lim_{k\to\infty}e^{\mu_k}=e^{\bh^{w_i}(m/n)}$ as required.
\end{proof}

\begin{example}
Note that entropy bounds obtained in this way depend only on local
features of the code of the periodic orbit under consideration.

For example, let~$R$ be any horseshoe periodic orbit for which
$\overline{c_R}$ contains the word $001\opti111\opti100$. Then some
shift of $\overline{c_R}$ is of the form $bv\cdot f$, where $v=\opti
111\opti = \opti w_1\opti$, and $q(b)\le 1/3$, $q(f)\le 1/3$. Thus
$r^{w_1}(R)\le 1/3$, and hence $h(R)\ge \bh^{w_1}(1/3)$, which is the
logarithm of the unique root in $x>1$ of the polynomial
$\overline{H}^1_{1/3}(x)$. 

Now $\overline{H}^1_{1/3}(x)=(x^3-1)g_1(x)+2(x^2-1)(x^6+1)$ (note that
$f_{1/3}(x)=0$), which simplifies (using $g_1(x)=x^5(x^3-x^2-x-1)-2$)
to
\[\overline{H}^1_{1/3}(x)={x}^{2} \left( x^4-1 \right)  \left( {x}^{5}-{x}^{4}-{x}^{3}+2\,x-2 \right) .\]

Thus any periodic orbit~$R$ whose code contains this word has
$h(R)\ge\log(1.47669)$. 

Compare this to the topological entropy of $P^{w_1}_{1/3}$ itself,
which is given by the largest positive root of the polynomial
\begin{eqnarray*}
H^1_{1/3}(x)&=&x^4g_1(x)-x^8g_1(1/x)\\
&=&\left( x^4-1 \right) \left( {x}^{8}-{x}^{7}-{x}^{6}-{x}^{5}+3\,{x}^{4}-{x}^{3}-{x}^{2}-x+1
 \right) 
,\end{eqnarray*}
giving $h^{w_1}_{1/3}\simeq \log(1.56294)$.
\end{example}

\bibliographystyle{amsalpha}

\bibliography{decbib}

\end{document}